\documentclass [a4paper,10pt]{article}
\usepackage [utf8] {inputenc}
\usepackage {amsmath}
\usepackage {amssymb}
\usepackage {amsthm}
\usepackage {amsfonts}
\usepackage {color}
\newtheorem {lemma}{Lemma}

\newtheorem {theorem}{Theorem}
\newtheorem* {theorem*}{Theorem}
\newtheorem {notation}{Notation}

\newcommand {\arrayfi }[1]{&&\hspace {20mm}{#1}\hspace {4mm}}
\title {A Montgomery-Hooley theorem for the $k$-fold divisor function}
\author {Tomos Parry}
\date {}
\begin {document}
\maketitle
\begin {abstract}
%Let $d_k(n)$ denote the $k$-fold divisor function.  %For all $k$, an asymptotic formula for 
%\[ V(x,Q):=\sum _{q\leq Q}\sum _{a=1}^q\left |\sum _{n\leq x\atop {n\equiv a(q)}}d(n)-\text { main term }\right |^2\]
%is established for a wide range of $q$.
Let $d_k(n)$ denote the $k$-fold divisor function.  For a wide range of large $q$ the expected bound
\[ \sum _{n\leq x\atop {n\equiv a(q)}}d_k(n)-\text { main term }\approx \sqrt {\frac {x}{q}}\]
is shown to be true in an average sense - for all $k$.  This generalises the work of Pongsriiam and Vaughan \cite {vaughandap} who studied $k=2$, and answers the work of Rodgers and Soundararajan \cite {sound}, who used the asymptotic large sieve to study a smoothed version of the problem. We use a circle method approach as developed by Goldston and Vaughan \cite {goldstonvaughan} to study the unsmoothed problem.
\end {abstract}
%We want to evaluate
%\begin {eqnarray*}
%S_t(Z)&=&\sum _{ql\leq y\atop {l<t/Z}}\frac {1}{l}\sum _{m<t-Zl\atop {t<n\atop {n-m=ql}}}nd(n)d(m)
%\\ &=&\int _0^1F(-\alpha )f^*(\alpha )f(-\alpha )d\alpha 
%\end {eqnarray*}
%where
%\[ F(\alpha )=\sum _{ql\leq y\atop {l<t/Z}}\frac {e(\alpha ql)}{l}\]
%\[ f^*(\alpha )=\sum _{t<n\leq y}nd(n)e(n\alpha )\]
%\[ f(\alpha )=\sum _{m<t-Zl}d(m)e(m\alpha ).\]
%The circle method analysis should give a main term
%\begin {eqnarray*}
%&&\sum _{q=1}^\infty G_q(q)^2\sum _{uv\leq y\atop {u<t/Z}}\frac {c_q(uv)}{u}\sum _{m<t-Zl\atop {t<n\atop {n-m=uv}}}\frac {d}{dn}\left \{ n^2(\log n+c)\right \} \frac {d}{dm}\left \{ m(\log m+c)\right \} 
%\end {eqnarray*}
%The sum is
%\begin {eqnarray*}
%&&2\sum _{m<t-Zl\atop {t<n\atop {n-m=uv}}}n\left (\log n+c+1/2\right )\left (\log m+c+1\right )
%\\ &&\hspace {10mm}=\hspace {4mm}\sum _{t,uv<n<y,t-Zu+uv}n\left (\log n+c+1/2\right )\left (\log (n-uv)+c+1\right )
%\\ &&\hspace {10mm}=:\hspace {4mm}\sum _{t,uv<n<y,t-Zu+uv}nf(n)f(n-uv)
%\\ &&\hspace {10mm}=:\hspace {4mm}\sum _{t,uv<n<y,t-Zu+uv}f_{uv}(n)
%\end {eqnarray*}
%so we have a total main term 
%\begin {eqnarray*}
%&&\sum _{q=1}^\infty G_q(q)^2\sum _{uv\leq y\atop {u<y/Z}}\frac {c_q(uv)}{u}\sum _{uZ,uv<n<y,y-Zu+uv}f_{uv}(n)\int _{Zu,n+Zu-uv}^{y,n}dt
%\\ &&\hspace {10mm}=\hspace {4mm}\sum _{q=1}^\infty G_q(q)^2\sum _{u<y/Z}c_q(uv)\sum _{Z<v\leq y/u}(v-Z)\sum _{uv<n<y}f_{uv}(n).
%\end {eqnarray*}
%From a lemma (Goldston variance document) the terms $v<Z$ can be added in with total error $\ll y^3Z$ and then the 
\begin {center}
\section {-\hspace {5mm}Introduction}
\end {center}
For $k\geq 1$ define the $k$-fold divisor function
%\[ d(n)=\sum _{d|n}1\]
%and more generally, for $k\geq 1$, the $k$-fold divisor function 
\[ d_k(n)=\sum _{u_1\cdot \cdot \cdot u_k=n\atop u_1,...,u_k\in \mathbb N}1.\]
The question of the distribution of values of $d_k(n)$ has a rich history, the values being linked to the Riemann zeta function through
\[ \zeta (s)^k=\sum _{n=1}^\infty \frac {d_k(n)}{n^s}.\]
%here, of course, the Riemann zeta function $\zeta (s)$, is given by
%\[ \zeta (s)=\sum _{n=1}^\infty \frac {1}{n^s},\]
%pivotal in analytic number theory due to the information it carries, with its Euler product, about the primes.
So too does the question of their distribution in arithmetic progressions, being in turn linked to Dirichlet $L$-functions and to, in particular, moments of Dirichlet $L$-functions, as indeed is said, for example, on page 3 of \cite {sound}.  Consequently larger $k$ often poses the more problems.  The driving question seems to have been \emph {the divisor problem for arithmetic progressions}, which asks when do we have a formula
\begin {eqnarray}\label {divisorproblemap}
\sum _{n\leq x\atop {n\equiv a(q)}}d_k(n)\approx \frac {xf_x(q,a)}{q}
\end {eqnarray}
or in other words when are the values of $d_k(n)$ uniformly distributed across the residue classes; here $f_x(q,a)$ is the residue at $s=1$ of something that looks like $\zeta (s)^kx^{s-1}/s$ and which we can think of therefore as having size no bigger than $(\log x)^{k-1}$.  %\[ \sum _{d|q}\frac {\mu (q/d)}{\phi (q/d)d^s}\sum _{n=1\atop {(n,q/d)=1}}^\infty \frac {d_k(dn)}{n^s}\]
%$f_s(q)\zeta (s)^kx^{s-1}/s$ where $f_s(q)$ . 
%The driving question seems to have been \emph {the divisor problem for arithmetic progressions}, which asks when do we have, for $(a,q)=1$, a formula
%\begin {eqnarray}\label {divisorproblemap}
%\sum _{n\leq x\atop {n\equiv a(q)}}d_k(n)\approx \frac {xf_x(q)}{\phi (q)}
%\end {eqnarray}
%or in other words when are the values of $d_k(n)$ uniformly distributed across the reduced residue classes; here $f_x(q)$ is the residue at $s=1$ of $L^k(s,\chi _0)x^{s-1}/s$. 
 % and as such is a polynomial in $\log x$ of degree $\leq k-1$ with coefficients $\alpha _q$ essentially bounded.  For a precise formulation see, for example, \cite {fi}.  
If \eqref {divisorproblemap} holds for all $q$ up to essentially $x^{\theta }$ then $d_k(n)$ is said to have \emph {exponent of distribution $\theta $}.  A critical value is $\theta =1/2$ which has direct consequences for the prime numbers - see Theorem 4 of \cite {fouvry} or, of course, consider the breakthrough ideas of Zhang on bounded gaps.  For $k=2$ it is by now a classical result of Selberg and of Hooley that $\theta =2/3$ is valid but further improvements have proven themselves elusive.  For $k=3$ important power-improvements on $\theta =1/2$ were made in the 80s by Friedlander and Iwaniec and for $k=4$ we have no better than $\theta =1/2$ from the 60s due to Linnik - see pages 32-33 of \cite {nguyen2} for references to these results.  For larger $k$ it is not known whether $\theta \geq 1/2$.  
\\
\\ So as is usual in these circles of questions one is inclined to ask instead what happens if we average - either by considering the variance
\[ V_x(q):=\sum _{a=1}^q\left |\sum _{n\leq x\atop {n\equiv a(q)}}d_k(n)-\frac {xf_x(q,a)}{q}\right |^2\]
%a result of the form $V_x(q)\approx x$ saying essentially the error in \eqref {divisorproblemap} is $\sqrt {x/q}$ on average, 
or the averaged variance
\[ V(x,Q):=\sum _{q\leq Q}V_x(q)\]
the point being that asymptotics for these of sizes around $x$ and $xQ$ respectively would confirm an error around $\sqrt {x/q}$ on average.  Starting with a result of Motohashi in the 70s there was a line of work (of which we mention \cite {banks}, \cite {blomer}, \cite {prapanpong}) on these variances culminating in 2012 with a very satisfactory result of Lau and Zhao \cite {lauzhao} which provided an asymptotic for $V_x(q)$ in a range essentially $\sqrt x\leq q\leq x$.  For $k>2$ however, an asymptotic formula for either variance isn't known, and in \cite {sound} Rodgers and Soundararajan put forward conjectures as to the true sizes of these variances for general $k$, these based on comparison with the function field case, investigated in \cite {keating}.  They then go on to provide evidence for these conjectures by approaching the problem with smooth weights.  Allowing smooth weights usually simplifies things (consider the usual Dirichlet divisor problem, where improvements on an error $\ll x^{1/3}$ are considered deep whilst with smoothing an error $\ll 1$ is straightfoward) and indeed they achieve for general $k$ an asymptotic formula for smooth $V(x,Q)$ in a range essentially $Q\geq x^{k/(k+2)}$.  More recently still, Nguyen \cite {nguyen1} has found an asymptotic formula even for smooth $V_x(q)$ in a range essentially $x^{1/(k-1)}\leq q\leq x^{1/k}$.  A step in a different direction is to change the approximating main term to a more workable but less natural or probabilistic one - for an asymptotic formula in this case see \cite {bf}.
\\
\\ But as the authors in \cite {sound} allude to themselves, it is another challenge to prove such results with sharp cut-offs, their concern in the paper being confirming the nature of the main term in the integer case at all.  But very recently progress has been done even in the sharp cut-off case, with Nguyen \cite {nguyen2} finding an asymptotic formula for $V(x,x)$ in the case $k=3$.  In this article, we take things further by providing one for $V(x,Q)$ for all $k$.    
\\
\\ Let's state our theorem.  In Lemma \ref {kfold} (C) we will prove that there is a quantity $f_x(q,a)$ such that for fixed $q,a\in \mathbb N$
\begin {equation}\label {ff}
\sum _{n\leq x\atop {n\equiv a(q)}}d_k(n)\sim \frac {xf_x(q,a)}{q}\hspace {10mm}x\rightarrow \infty .
\end {equation}
(For a precise definition see that lemma, but essentially $f_x(q,a)$ is a degree $\leq k-1$ polynomial in $\log x$ with coefficients $\ll _\epsilon q^\epsilon $).  Suppose we know that for some $\mathfrak c\in (1/2,1)$ we have 
%\[ K(s)\ll \frac {1}{|s|^2|s+1|}.\]
\[ \int _{\pm \infty }\frac {\max _{0\leq n\leq 2k-2}|\zeta ^{(n)}(\sigma +it)|^{k(k-2)}dt}{(1+|t|)^{3/2+\mathfrak c}}\ll 1\hspace {10mm}\text {for }\sigma \geq \mathfrak c.\]
A bound essentially $\zeta (s)\ll t^{(1-\sigma )/2}$ is classical, and therefore from Cauchy's integral formula so is one $\zeta ^{(n)}(s)^k\ll t^{k(1-\sigma )/2}$, so we can certainly take $\mathfrak c\in (1-1/k(k-2),1)$ for $k>2$ or take $\mathfrak c\in (1/2,1)$ for $k=2$.  
Let $\mathfrak d\in (0,1)$ be any value for which we know
\begin {equation}\label {square}
\sum _{n\leq x}d_k(n)^2=xP(\log x)+\mathcal O_\epsilon \left (x^{1-\mathfrak d+\epsilon }\right )
\end {equation}
for some polynomial $P$ of degree $k^2-1$.  By page 12 of \cite {sound} we have for $\sigma >1$
\[ \sum _{n=1}^\infty \frac {d_k(n)^2}{n^s}=\zeta (s)^{k^2}f(s)\]
where $f$ for $\sigma \geq 1/2+\epsilon $ is holomorphic and $\ll _\epsilon 1$, so a Perron's formula argument shows $\mathfrak d=2/(k^2+4)$ to be permissable for general $k$ whilst a classical result of Wilson and Ramanujan says $\mathfrak d=1/2$ is okay for $k=2$.
\\
\\ Then our conclusion is
\\ 
\begin {theorem}\label {mh}
%Let $k\geq 2$ be in $\mathbb N$.    
Let $f_x(q,a)$ be as in \eqref {ff} and define
\begin {eqnarray*}
%x\in \mathbb R,\hspace {2mm}q,d\in \mathbb N\hspace {10mm}\mathcal M_x(q,d)&=&Res_{s=1}\Bigg \{ \frac {x^{s-1}}{s}\sum _{n=1\atop {(n,q)=d}}^\infty \frac {d_k(n)}{n^{s}}\Bigg \} 
%\\ x\in \mathbb R,\hspace {2mm}q,a\in \mathbb N\hspace {10mm}
E_x(q,a)=\sum _{n\leq x\atop {n\equiv a(q)}}d_k(n)-\frac {xf_x(q,a)}{q}\hspace {6mm}\text { and }\hspace {6mm}V(x,Q)&=&\sum _{q\leq Q}\sum _{a=1}^q|E_x(q,a)|^2.
%We will prove in Lemma \ref {kfold} (B) that there is a quantity $\mathcal M_X(q,a)$ such that, for fixed $k\geq 2$ in $\mathbb N$ and fixed $q,a\in \mathbb N$, 
%\[ \sum _{n\leq x\atop {n\equiv a(q)}}d_k(n)\sim \frac {x\mathcal M_x(q,a)}{q}\hspace {10mm}\text {as }x\rightarrow \infty .\]
%\\ x,Q\in \mathbb R\hspace {10mm}V(x,Q)&=&\sum _{q\leq Q}\sum _{a=1}^q|E_x(q,a)|^2.
\end {eqnarray*}
Then there is a polynomial $P(\cdot ,\cdot )$ of degree $\leq k^2-1$ such that for $1\leq Q=o(x)$
\begin {eqnarray*}
V(x,Q)=xQP(\log x,\log Q)+\mathcal O_{k,\epsilon }\left (Q^2\left (\frac {x}{Q}\right )^{\mathfrak c}+\underbrace {x^{3/2+\epsilon }}_{k=2}+\underbrace {x^{2-4/(6k-3)+\epsilon }}_{k>2}+Qx^{1-\mathfrak d+\epsilon }\right ).
\end {eqnarray*}
%For $k=2$ the second error may be replaced by $x^{3/2+\epsilon }$.
\end {theorem}
%We end the introduction by comparing this result for $d_k$ against the same result for other sequences and for $d_2$, and by briefly discussing our proof.  The study of the mean-square variance averaging over both modulus and residue started with the Barban-Davenport-Halberstam Theorem, which provided a bound of the right size for the variance of the primes.  This was refined to an asymptotic formula both by Montgomery and by Hooley and correspondingly results of our kind are sometimes referred to as Montgomery-Hooley theorems.  Hooley, in a long series of papers of which we mention \cite {hooley}, and Vaughan \cite {vaughangeneral} inquired as to sequences aside from the primes, the most notable perhaps being the values of the M\" obius function and the squarefree numbers.   Initially the main ingredient was the large sieve but Goldston and Vaughan developed a new approach in \cite {goldstonvaughan}, which we follow here.   
\hspace {1mm}
\\ We consider the main feature of our theorem to be its validity for general $k$.  
\\
\\ With some work the polynomial could be explicitly worked out of course.  We note that we have chosen to average over a complete residue set, and not over a reduced one, as in most of the works referenced above, in agreement with a remark of Motohashi on page 178 of \cite {motohashi}.
\\
\\ The study of the mean-square of a function averaging over both residue and modulus started in the 60s with the Barban-Davenport-Halberstam Theorem, which provided an upper bound of the right size for the variance of the primes.  This was refined to an asymptotic formula by Montgomery and by Hooley and correspondingly results of our kind are sometimes referred to as Montgomery-Hooley theorems - see \cite {survey} for a survey.  In \cite {goldstonvaughan} Goldston and Vaughan developed a circle method approach to establishing Montgomery-Hooley theorems and it is that which we follow here.
\\
\\ Sequences or functions for which there exist Montgomery-Hooley theorems seem to have so far been mostly assumed to satisfy 
\begin {equation}\label {hwyr}
\sum _{n\leq x\atop {n\equiv a(q)}}F(n)\approx \frac {xf(q,a)}{q}
\end {equation}
although recently thinner sequences have also been investigated - see \cite {bruedernvaughan}.  The divisor function poses a different problem in that, as said on page 86 of \cite {survey}, the RHS of \eqref {divisorproblemap} doesn't factorise in the simple way of \eqref {hwyr}.  With this paper we are showing that the method of Goldston and Vaughan still succeeds.  
\\
\\ Finally we remark that as far as only upper bounds are concerned results will already be known even for $d_k(n)$, even if nothing explicit is stated in the literature - a general principle is that, through the large sieve, a Siegel-Walfisz theorem gives a result for the variance.  So just for the record we also formulate
\begin {theorem}\label {largesieve}
Let $V(x,Q)$ be as in Theorem \ref {mh}.  For $1\leq Q\leq x$
%\[ \sum _{q\leq Q}\sideset {}{'}\sum _{a=1}^q|E_X(q,a)|^2\ll x^\epsilon \left (xQ+\underbrace {x^{2-1/k}}_{k>2}\right ).\]
\[ V(x,Q)\ll _\epsilon x^\epsilon \left (xQ+\underbrace {x^{2-1/k}}_{k>2}\right ).\]
\end {theorem}
\hspace {1mm}
\\ and prove it following Lemma \ref {barbandavenporthalberstam}.  The proof of Theorem \ref {mh} starts on page 15.
%Hooley, in a long series of papers of which we mention \cite {hooley}, and Vaughan \cite {vaughangeneral} inquired as to sequences aside from the primes, the most notable perhaps being the values of the M\" obius function and the squarefree numbers.  Initially the main ingredient was the large sieve but Goldston and Vaughan developed a new approach in \cite {goldstonvaughan}, which we follow here.  
\begin {center}
\section {-\hspace {5mm}Lemmas}
\end {center}
We start with simple estimates for the $k$-fold divisor function's exponential sum and counts in arithmetic progressions.  We make use of a result of \cite {vaughandap} but otherwise our estimates are very simple - we are interested only in any power-savings.  In particular any improvements here would also work their way through to our theorems.
\\ 
\begin {lemma}\label {kfold}
Let $k,q,a,d\in \mathbb N$ with $k\geq 2$ and $d|q$ and let $X\geq 1$.  All error terms below are $\mathcal O_{k,\epsilon }((qX)^\epsilon \cdot \cdot \hspace {1mm}\cdot )$ and we write
\[ \theta =2/(k+1)\hspace {7mm}\Delta =1/(k-1)\hspace {7mm}\delta =2/(2k-1).\]
Then
\begin {alignat*}{2}
&(A)\hspace {23mm}\sum _{n\leq X\atop {(n,q)=d}}\frac {d_k(n)}{n}&\hspace {1mm}=\hspace {1mm}&Res_{s=0}\Bigg \{ \frac {X^{s}}{s}\sum _{n=1\atop {(n,q)=d}}^\infty \frac {d_k(n)}{n^{s+1}}\Bigg \}+\mathcal O\left (\frac {1}{X}\left (\frac {X}{d}\right )^{1-\theta }\right )%\mathcal O\left (\frac {1}{X}\left (\left (\frac {X}{d}\right )^{1-\theta }+1\right )\right )
\\ &(B)\hspace {15mm}\sum _{n\leq X}d_{k}(n)e\left (\frac {na}{q}\right )&\hspace {1mm}=\hspace {1mm}&\frac {X\mathcal M_X\left (q/(q,a)\right )}{q/(q,a)}+\mathcal O\left (X\left (\left (\frac {q}{X}\right )^{\Delta }+\frac {1}{\sqrt q}\left (\frac {q}{X}\right )^{\Delta /2}\right )\right )
\end {alignat*}
where 
\begin {eqnarray*}
%\\ g_X(q)&=&P(x,q)\hspace {5mm}\text { if }k=2
%\\ g_X(q)&=&\sum _{h_1,...,h_{k-2}|q}\frac {\mu (h_1)\cdot \cdot \cdot \mu (h_{k-2})P(X,q,\mathbf r,\mathbf h)}{h_1\cdot \cdot \cdot h_{k-2}}\hspace {5mm}\text { if }k>2
\mathcal M_X(q)&=&P(X,q)\hspace {5mm}\text {for }k=2
\\ \mathcal M_X(q)&=&\sum _{r_1|\cdot \cdot \cdot |r_{k-2}|q\atop {h_1|r_1,...,h_{k-2}|r_{k-2}}}\frac {\mu (h_1)\cdot \cdot \cdot \mu (h_{k-2})P(X,q,\mathbf r,\mathbf h)}{h_1\cdot \cdot \cdot h_{k-2}}\hspace {5mm}\text {for }k>2
%P(x,q)\hspace {5mm}\text { if }k=2
%\\ f_X(q)&=&\mathcal N_X^{}(q,q)
%\sum _{r_1|\cdot \cdot \cdot |r_{k-2}|q\atop {h_1|r_1,...,h_{k-2}|r_{k-2}}}\frac {\mu (h_1)\cdot \cdot \cdot \mu (h_{k-2})P(X,q,\mathbf r,\mathbf h)}{h_1\cdot \cdot \cdot h_{k-2}}
%\hspace {5mm}\text { if }k>2
%\\ \mathcal M_X(q,a)&=&\sum _{d|q}\frac {c_d(a)g_X(d)}{d}
%where $c_q(a)$ is Ramanujan's sum.  
\end {eqnarray*}  % and write
%\[ \mathfrak E_k(X)=\left \{ \begin {array}{ll}\sqrt X&\text { if }k=2%\\ \sqrt {Xq}&\text { if }k=3
%\\ X^{1-\Delta }q^\Delta &\text { if }k\geq 3.\end {array}\right \} %\hspace {5mm}\text{and}\hspace {5mm}\mathfrak F_k(X)=\left \{ \begin {array}{ll}X^{1/3}+\sqrt q&\text { if }k=2%\\ \sqrt {Xq}&\text { if }k=3
%\\ &\text { if }k=3
%\\ X^{1-\Delta }q^\Delta &\text { if }k\geq 4\end {array}\right \} 
%\]
%with $\Delta =\min \{ 1/2,1/(k-2)\} $.  
for some polynomials $P(\cdot ,...,\cdot )$ in $\log (\cdot ),...,\log (\cdot )$ of degree $\leq k-1$.
\\
\\ Let $\mathcal M_X(q)$ be as in claim (B), write
\[ c_q(n)%=\sideset {}{'}\sum _{a=1}^qe\left (\frac {an}{q}\right )
=\sum _{d|q,n}\mu (q/d)d=\sum _{a=1\atop {(a,q)=1}}^qe\left (\frac {an}{q}\right )\hspace {7mm}\text {for Ramanujan's sum, and let}\hspace {7mm}f_X(q,a)=\sum _{d|q}\frac {c_{d}(a)\mathcal M_X(d)}{d}.\]
Then
\begin {alignat*}{2}&(C)\hspace {25mm}\sum _{n\leq X\atop {n\equiv a(q)}}d_{k}(n)&\hspace {1mm}=\hspace {1mm}&\frac {Xf_X(q,a)}{q}+\mathcal O\left (X\left (\left (\frac {\sqrt q}{X}\right )^{\Delta }+\frac {1}{X^{\delta }}\right )\right )
\\ &(D)\hspace {19mm}\sum _{n\leq X}d_{k}(n)c_q(n)&\hspace {1mm}=\hspace {1mm}&\frac {X\phi (q)\mathcal M_X(q)}{q}+\mathcal O\left (X^{1-\theta }q^\theta \right ).
\end {alignat*}
\end {lemma}
\begin {proof}
Throughout we drop $X^\epsilon ,q^\epsilon $ factors from error terms, write $c=1/\log (2X)$, write $T$ for any parameter $1\leq T\ll X$, write $\Sigma _{\chi (q)}$ for a sum over the Dirichlet characters of modulus $q$ and write $L_\chi (s)$ for the Dirichlet $L$-function of a character $\chi $.  Before we address the individual claims of the lemma let us do some preparatory work in establishing \eqref {residue2} and \eqref {gcd2}-\eqref {wy3} below.  For $j\in \mathbb N$ write $P_j(\cdot ,...,\cdot )$ to mean simply any polynomial in $\log (\cdot ),...,\log (\cdot )$ of degree $\leq j$ and for any vector $\mathbf v$ write $\mathcal N_\mathbf v^{k}(q,d)$ for any quantity of the form
\[ %\mathcal N_{\mathbf v}^K(q,d)=
\sum _{r_kh_k|q\atop {...\atop {r_1\cdot \cdot \cdot r_kh_k|q\atop {r_1\cdot \cdot \cdot r_k|d\atop {H|q/d}}}}}\frac {\mu (h_1)\cdot \cdot \cdot \mu (h_k)\mu (H)P_{k+1}(\mathbf v,q,d,\mathbf r,\mathbf h,H)}{h_1\cdot \cdot \cdot h_kH}\]
and write $\mathcal N_\mathbf v^0(q,d)$ for one of the form
\[ \sum _{H|q/d}\frac {\mu (H)P_{1}(\mathbf v,q,d,H)}{H}.\]
%so that in particular
%\begin {equation}\label {nn}
%\mathcal M_X(q)=\mathcal N_X^{k-2}(q,q).%\hspace {8mm}\text {and}\hspace {8mm}\frac {\phi (q)}{q}\mathcal N_\mathbf v^{k-1,k}(q,1)=\mathcal N_\mathbf v^{k,k}(q,1).
%\end {equation}
It is straightforward to establish
\[ \sum _{n=1\atop {(n,q)=d}}^\infty \frac {d_2(n)}{n^s}=\frac {\zeta (s)^2}{d^s}\sum _{rh|q\atop {r|d\atop {H|q/d}}}\frac {\mu (h)\mu (H)}{(hH)^s}\]
and then induction shows
\[ \sum _{n=1\atop {(n,q)=d}}^\infty \frac {d_k(n)}{n^{s}}=\frac {\zeta (s)^{k}}{d^s}\sum _{r_kh_k|q\atop {\cdot \cdot \cdot \atop {r_1\cdot \cdot \cdot r_{k-1}h_{k-1}|q\atop {r_1\cdot \cdot \cdot r_{k-1}|d\atop {H|q/d}}}}}\frac {\mu (h_1)\cdot \cdot \cdot \mu (h_{k-1})\mu (H)}{(h_1\cdot \cdot \cdot h_{k-1}H)^s}\]
so
%\begin {eqnarray}\label {residue1}
%\[ Res_{s=1}\Bigg \{ \frac {X^{s-1}}{s}\sum _{n=1\atop {(n,q)=d}}^\infty \frac {d_k(n)}{n^s}\Bigg \} =\frac {\mathcal N_X^{k-1,k-1}(q,d)}{d}\hspace {10mm}\]
%\end {eqnarray}
\begin {equation}\label {residuepoly}
Res_{s=1}\Bigg \{ \frac {X^{s-1}}{s}\sum _{n=1\atop {(n,q)=d}}^\infty \frac {d_k(n)}{n^s}\Bigg \}  \hspace {5mm}\text { is a polynomial in $\log X$ of degree $\leq k-1$} 
\end {equation}
and
\begin {eqnarray}\label {residue2}
Res_{s=0}\Bigg \{ \frac {X^{s}}{s}\sum _{n=1\atop {(n,q)=d}}^\infty \frac {d_k(n)}{n^{s+1}}\Bigg \} =\frac {\mathcal N_X^{k-1}(q,d)}{d}.\hspace {10mm}
\end {eqnarray}
Orthogonality and Perron's formula says that for $(q,a)=1$ 
%and 
%We will make use of $\phi /id=\mu /id\star 1$ and of
%\[ \sum _{d|n}\mu (d)=\left \{ \begin {array}{ll}1&\text { if }n=1\\ 0&\text { if }n>1\end {array}\right .\]
%which means a coprimality condition may vanish from one line to the next. 
%\begin {eqnarray}\label {wy}
%&&\sum _{v<X/yq}\left \{ \frac {Xf_{X/v}^{k,k-1}\left (q/(q,v)\right )}{v}-yf_y^{k-1}\left (q/(q,v)\right )\right \} \notag 
%\\ &&\sum _{v<XR/yq\atop {(v,R)=1}}\left \{ \frac {Xf_{XR/vq}^{k,k-1}(R)}{vq/d}-yf_y^{k-1}(R)\right \} \notag 
%\\ \arrayfi =\sum _{\mathbf h|\mathbf r|d}\frac {\mu (\mathbf h)}{\mathbf h}\sum _{H|R}\frac {\mu (H)}{H}\sum _{v<XR/yqH}\left \{ \frac {XP_{k-1}(X,v,q,R,H,\mathbf r,\mathbf h)}{vq/R}-yP_{k-1}(y,q,\mathbf r,\mathbf h)\right \} \notag 
%\\ \arrayfi =\frac {X}{q/R}\sum _{r_1|\cdot \cdot \cdot |r_{k-2}|R\atop {h_1|r_1,...,h_{k-2}|r_{k-2}}}\frac {\mu (\mathbf h)}{\mathbf h}\sum _{H|R}\frac {\mu (H)P_{k}(X,y,q,\mathbf r,\mathbf h,H,R)}{H}+\mathcal O\left (y\right )
%\hspace {10mm}
%\end {eqnarray}
\begin {eqnarray*}
\sum _{n\leq X\atop {n\equiv a(q)}}d_k(n)&=&\frac {1}{\phi (q)}\sum _{\chi (q)}\overline \chi (a)\int _{1+c\pm iT}\frac {L_\chi (s)^kX^sds}{s}+\mathcal O\left (\frac {X}{T}\right )
%\\ &=&\frac {1}{\phi (q)}Res_{s=1}\Bigg \{ \frac {(X/d)^s}{s}\sum _{n=1\atop {(n,q)=1}}^\infty \frac {d_k(n)}{n^{s}}\Bigg \} +\mathcal O\left (qX^{1-2/k}\right ).
\end {eqnarray*}
so, using the well-known bound $L_\chi (s)\ll _qt^{(1-\sigma )/2}$ to move the integral to the left at the cost of errors
\begin {eqnarray*}
\left (\int _{\epsilon \pm iT}+\int _{\epsilon +iT}^{c+iT}\right )\frac {L_\chi (s)^kX^sds}{s}&\ll _q&T^{k/2}+\frac {X}{T}\ll X^{1-2/(k+2)}
\end {eqnarray*}
and picking up a residue we may say
\begin {eqnarray}\label {coprime}
\text {for }(q,a)=1\hspace {10mm}\sum _{n\leq X\atop {n\equiv a(q)}}d_k(n)&=&\frac {X}{\phi (q)}Res_{s=1}\Bigg \{ \frac {L_{\chi _0}(s)^kX^{s-1}}{s}\Bigg \} +\mathcal O_q\left (X^{1-2/(k+2)}\right ).\hspace {10mm}
%\\ &=&\frac {1}{\phi (q)}Res_{s=1}\Bigg \{ \frac {(X/d)^s}{s}\sum _{n=1\atop {(n,q)=1}}^\infty \frac {d_k(n)}{n^{s}}\Bigg \} +\mathcal O\left (qX^{1-2/k}\right ).
\end {eqnarray}
To remove the coprimality condition we use an argument from \cite {hb}.  Define Popovici's function $\mu ^k=\mu \star \cdot \cdot \cdot \star \mu $ with $k$ convolution factors and let $\delta $ be a power of a prime $p$.  We have $N>k\implies \mu ^k(p^N)=0$ and $N\geq 2\implies (d_2(\delta \cdot )\star \mu ^2)(p^N)=0$.  Define $c_\delta =d_{k+1}(\delta /p)$ so that on prime powers $d_{k+1}(\delta \cdot )=d_k(\delta \cdot )\star 1+c_\delta 1$ and therefore $d_{k+1}(\delta \cdot )\star \mu ^{k+1}=d_k(\delta \cdot )\star \mu ^k+c_\delta \mu ^k$ so with the last sentence 
\begin {equation}\label {dolig}
N\geq k\implies (d_k(p^D\cdot )\star \mu ^k)(p^N)=0
\end {equation}
holds for all $k\in \mathbb N$, all primes $p$, and all $N,D\geq 1$.  For $\delta \in \mathbb N$ define $F_\delta =d_k(\delta \cdot )\star \mu ^k$ which is multiplicative in the sense that if $D_i,N_i$ are non-negative powers of primes $p_i$ then
\[ F_{D_1\cdot \cdot \cdot D_r}(N_1\cdot \cdot \cdot N_r)=F_{D_1}(N_1)\cdot \cdot \cdot F_{D_r}(N_r)\]
%$(n\delta ,n'\delta ')\implies F_{\delta \delta '}(nn')=F_\delta (n)F_{\delta '}(n)$ 
so \eqref {dolig} says $F_\delta (n)=0$ unless $n|\delta ^{k-1}$, and therefore
%\begin {eqnarray}\label {gcd}
\[ \sum _{n}d_k(\delta n)=\sum _{n,m\atop {m|\delta ^{k-1}}}F_\delta (m)d_k(n)\hspace {5mm}\text { and }\hspace {5mm}\sum _{m>X}\frac {F_\delta (m)}{m}\ll \frac {\delta ^\epsilon }{X}.\]
%\end {eqnarray}
Consequently if for coprime $q,a\in \mathbb N$ and $Q\in \mathbb N$ with $(Q,q)=1$ we have for some $0\leq w,E\leq 1$ 
\[ \sum _{n\leq X\atop {n\equiv a(q)\atop {(n,Q)=1}}}\frac {d_k(n)}{n^w}=\frac {1}{\phi (q)}Res_{s=1-w}\Bigg \{ \frac {X^{s}}{s}\sum _{n=1\atop {(n,qQ)=1}}^\infty \frac {d_k(n)}{n^{s+w}}\Bigg \} +\mathcal O\left (X^{1-w-E}\right )\]
then for general $q,a,Q,D\in \mathbb N$ with $D|Q$ and $(Q,q)=1$ we have, writing $\delta =(q,a)\hspace {2mm}q=\delta Dq'\hspace {2mm}a=\delta Da'$,% and assuming $(\delta ,Q)=1$,
\begin {eqnarray}\label {gcd2}
\sum _{n\leq X\atop {n\equiv a(q)\atop {(n,Q)=D}}}\frac {d_k(n)}{n^w}&=&%\sum _{n\leq x}\frac {1}{n}\sum _{m|n}F_q(m)d_k(n/m)
%\\ &=&\sum _{nm\leq x}\frac {F_q(m}d_k(n){nm}
%\\ &=&
\frac {1}{(\delta D)^w}\sum _{m\leq X/\delta D\atop {m|(\delta D)^{k-1}}}\frac {F_{\delta D}(m)}{m^w}\sum _{n\leq X/m\delta D\atop {nmD\equiv a'(q')\atop {(nm\delta ,Q/D)=1}}}\frac {d_k(n)}{n^w}\notag 
\\ &=&\frac {1}{(\delta D)^w\phi (q')}\sum _{m\leq X/\delta D\atop {m|(\delta D)^{k-1}}}\frac {F_{\delta D}(m)}{m^w}Res_{s=1-w}\Bigg \{ \frac {(X/m\delta D)^{s}}{s}\sum _{n=1\atop {(nmD,q')=1\atop {(nm\delta ,Q/D)=1}}}^\infty \frac {d_k(n)}{n^{s+w}}\Bigg \} \notag 
\\ &&\hspace {45mm}+\hspace {4mm}\mathcal O\left (\frac {1}{(\delta D)^w}\left (\frac {X}{\delta D}\right )^{1-w-E}\right )\notag 
\\ &=&\frac {1}{\phi (q')}Res_{s={1-w}}\Bigg \{ \frac {X^{s}}{s}\sum _{n=1\atop {(n,q)=\delta \atop {(n,Q)=D}}}^\infty \frac {d_k(n)}{n^{s+w}}\Bigg \} %Res_{s=1}\Bigg \{ \frac {X^{s-1}}{s}\sum _{n=1\atop {(n,q)=d}}^\infty \frac {d_k(n)}{n^{s}}\Bigg \} 
+\mathcal O\left (\frac {1}{(\delta D)^w}\left (\frac {X}{\delta D}\right )^{1-w-E}+\frac {1}{X^w}\right ).
\end {eqnarray}
In particular \eqref {coprime} and \eqref {residuepoly} give
\begin {eqnarray}\label {residue}
%\text {for a polynomial of degree $\leq k-1$}\hspace {5mm}
\sum _{n\leq X\atop {n\equiv a(q)}}d_k(n)=X\cdot \underbrace {\frac {1}{\phi \left (q'\right )}Res_{s=1}\Bigg \{ \frac {X^{s-1}}{s}\sum _{n=1\atop {(n,q)=\delta }}^\infty \frac {d_k(n)}{n^s}\Bigg \} }_{\text {polynomial in $\log X$ of degree $\leq k-1$ }}+\mathcal O_q\left (X^{1-2/(k+2)}\right ).\hspace {10mm}
%=\frac {1}{d}\sum _{r_1h_1|q\atop {...\atop {r_1\cdot \cdot \cdot r_{k-1}h_{k-1}|q\atop {r_1\cdot \cdot \cdot r_{k-1}|d\atop {H|q/d}}}}}\frac {\mu (\mathbf h)\mu (H)P^*(X/d,\mathbf h,H)}{\mathbf hH}
\end {eqnarray}
Finally the Euler-Maclaurin summation formula says that
%so that in particular
%\begin {equation}\label {residue}
%\mathcal M_X^k(q,d)=\mathcal N_X^{k-1}(q,d).
%\end {equation}
\begin {eqnarray*}
\sum _{v<X}\frac {P_j(v)}{v}%&=&\int _{1/2}^{X}\frac {P_j(t)dt}{t}+\int _{1/2}^\infty \frac {B_1(\{ t\} )P_j(t)dt}{t^2}%+\mathcal O\left (\frac {1}{Y}\right )
%\\ &=&\int _{ZY}^{ZX}\frac {P(\log U)dU}{U}
%\\ &=&
%\\ &=&
=P_{j+1}(X)+\mathcal O\left (\frac {1}{X}\right )
\end {eqnarray*}
%so 
%\begin {eqnarray*}
%\sum _{v<X\atop {(v,R)=1}}\frac {P_d(v)}{v}%&=&\sum _{H|R}\frac {\mu (H)}{H}\sum _{v<X/H}\frac {P_d(Hv)}{v}
%\\ &=&
%\\ &=&
%=\sum _{H|R}\frac {\mu (H)P_{d+1}(X,H)}{H}%+\mathcal O\left (\frac {1}{X}\right )
%\end {eqnarray*}
so for any $y>0$
\begin {eqnarray}\label {wy2}
\sum _{v<X/y}\frac {1}{q/(q,v)}\left \{ \frac {X}{v}\mathcal N^{k-2}_{X/v}\left (\frac {q}{(q,v)},\frac {q}{(q,v)}\right )-y\mathcal N^{k-2}_{y}\left (\frac {q}{(q,v)},\frac {q}{(q,v)}\right )\right \} %\notag 
%\\ &&\hspace {20mm}=
=\frac {X\mathcal N^{k-1}_{X,y}(q,q)}{q}+\mathcal O(y)%\hspace {10mm}
%\\ \arrayfi =\sum _{h_1,...,h_{k}|q}\frac {\mu (\mathbf h)}{\mathbf h}\sum _{H|q}\frac {\mu (H)}{H}\sum _{v<X/yH}\left \{ \frac {XP_k(X,v,q,d,\mathbf h,H)}{v}-yP_k(y,d,\mathbf h)\right \} \notag 
%\\ \arrayfi =X\sum _{h_1,...,h_{k}|q}\frac {\mu (\mathbf h)}{\mathbf h}\sum _{H|q}\frac {\mu (H)P_{k}(X,y,q,\mathbf h,d,H)}{H}+\mathcal O\left (y\right )
%\\ \arrayfi =X\mathcal N^{k+1}_{X}(q,1)
%\hspace {10mm}
\end {eqnarray}
and
\begin {eqnarray}\label {wy3}
&&\sum _{v<X/y\atop {(v,q)=1}}\left \{ \frac {X}{v}\mathcal N^{k-1,k-1}_{X/v}\left (q,1\right )-y\mathcal N^{k-1,k-1}_{y}\left (q,1\right )\right \} =X\mathcal N^{k,k}_{X,y}(q,1)+\mathcal O(y)\hspace {10mm}
%\\ \arrayfi =\sum _{h_1,...,h_{k}|q}\frac {\mu (\mathbf h)}{\mathbf h}\sum _{H|q}\frac {\mu (H)}{H}\sum _{v<X/yH}\left \{ \frac {XP_k(X,v,q,d,\mathbf h,H)}{v}-yP_k(y,d,\mathbf h)\right \} \notag 
%\\ \arrayfi =X\sum _{h_1,...,h_{k}|q}\frac {\mu (\mathbf h)}{\mathbf h}\sum _{H|q}\frac {\mu (H)P_{k}(X,y,q,\mathbf h,d,H)}{H}+\mathcal O\left (y\right )
%\\ \arrayfi =X\mathcal N^{k+1}_{X}(q,1)
%\hspace {10mm}
\end {eqnarray}
and now we are ready to turn to the individual claims.  \textbf {(A)} Let $Q\in \mathbb N$.  Perron's formula says that 
\begin {eqnarray}\label {perronsformula}
\sum _{n\leq X\atop {(n,Q)=1}}\frac {d_k(n)}{n}&=&\int _{c\pm iT}\frac {L_{\chi _0}(s+1)^kX^{s}ds}{s}+\mathcal O\left (\frac {1}{T}\right ).
\end {eqnarray} 
Arguing as in Section 12.2 of \cite {rzf} it may be shown that
\begin {eqnarray*}
\int _{-1-\epsilon \pm iT}\frac {L_{\chi _0}(s+1)^kX^{it}ds}{s}%&=&\sum _{n=1\atop {}}^\infty \frac {d_k(n)}{n^{1+\epsilon }}\int _{-1-\epsilon \pm iT}\frac {\chi ^k(s)(Xn)^{it}ds}{s}
\ll T^{k/2-1/2}
%\end {eqnarray*}
\hspace {10mm}\text {whilst easily}\hspace {10mm}%with the bound $\zeta (s+1)^k\ll ^{-k\sigma /2}$ for $\sigma \in [-1,c]$,
%\begin {eqnarray*}
\int _{-1-\epsilon +iT}^{c+iT}|L_{\chi _0}(s+1)|^kX^\sigma ds\ll \frac {T^{k/2}}{X}+1
\end {eqnarray*}
so moving the integral in \eqref {perronsformula} to the left we pick up a residue and introduce errors
\begin {eqnarray*}
\ll \frac {T^{k/2-1/2}}{X}+\frac {1}{T}\ll X^{-\theta }
\end {eqnarray*}
to conclude %for $X\geq 2$
\begin {eqnarray*}%\label {now}
\sum _{n\leq X\atop {(n,Q)=1}}\frac {d_k(n)}{n}=Res_{s=0}\Bigg \{ \frac {X^s}{s}\sum _{n=1\atop {(n,Q)=1}}^\infty \frac {d_k(n)}{n^{s+1}}\Bigg \} +\mathcal O\left (X^{-\theta }\right ).
\end {eqnarray*}
%we have% with \eqref {noq}
%\begin {eqnarray*}
%\sum _{n\leq X/d\atop {(n,q/d)=1}}\frac {d_k(dn)}{n}&=&%\sum _{n\leq x}\frac {1}{n}\sum _{m|n}F_q(m)d_k(n/m)
%\\ &=&\sum _{nm\leq x}\frac {F_q(m}d_k(n){nm}
%\\ &=&
%\sum _{m\leq X/d\atop {m|d^{k-1}\atop {(m,q/d)=1}}}\frac {F_d(m)}{m}\sum _{n\leq X/dm\atop {(n,q/d)=1}}\frac {d_k(n)}{n}
%\\ &=&\sum _{m|d^{k-1}\atop {(m,q/d)=1}}\frac {F_q(m)}{m}Res_{s=0}\Bigg \{ \frac {L_{\chi _0}(s+1)^k(X/dm)^s}{s}\Bigg \} +\mathcal O\left (\left (\frac {X}{d}\right )^{-\theta }+\frac {d}{X}\right )
%\\ &=&Res_{s=0}\Bigg \{ \frac {(X/d)^s}{s}\sum _{n=1\atop {(n,q)=d}}^\infty \frac {d_k(n)}{n^{1+s}}\Bigg \} +\mathcal %O\left (\left (\frac {X}{d}\right )^{-\theta }+\frac {d}{X}\right )
%\end {eqnarray*} 
%and we're done.  
From this, the sentence containing \eqref {gcd2}, and \eqref {residue2} we're done.  \textbf {(B)} Assume $(a,q)=1$ and $q\ll X$.  The claim in question is true for $k=2$ by Lemma 3.2 of \cite {prapanpong} so let's suppose it's true for some $k\geq 2$ and aim to show the claim for $k+1$; this means we want a main term of the form $\mathcal N_X^{k-1}(q,q)$.  Write $D^k_{a/q}(X)$ for the sum in question so that for a parameter $0<y\leq X$ we have
\begin {eqnarray}\label {split}
D^{k+1}_{a/q}(X)&=&\sum _{u\leq y}d_k(u)\sum _{v\leq X/u}e\left (\frac {uva}{q}\right )+\sum _{v<X/y}\Big (D_{va/q}^k(X/v)-D_{va/q}^k(y)\Big ).\hspace {10mm}
%\\ &=&\frac {X}{q}\sum _{r_1,...,r_{k-1}|q}\frac {\phi (\mathbf r)P(x/q,\mathbf r)}{\mathbf r}+\mathcal O\left (y+q\right )
%\\ &&\hspace {0mm}+\hspace {4mm}\frac {x}{q}\sum _{r_1|\cdot \cdot \cdot |r_{k-2}|R|q\atop {\mathbf h|\mathbf r\atop {H|R}}}\frac {\mu (\mathbf h)\mu (H)P_k(x,y,q,R,\mathbf r)}{\mathbf hH}+\mathcal O\left (x^{1/(1+\delta )}\right )
%\\ \arrayfi {=:}\frac {xH_y(q)}{q}+\mathcal V+\mathcal O\left (\frac {x}{y^\delta }+y+q\right ).
\end {eqnarray}
As
\begin {eqnarray*}
\sum _{u\leq y\atop {q\nmid u}}d_k(u)\sum _{v\leq x/u}e\left (\frac {uva}{q}\right )%&\ll &\sum _{u\leq y\atop {q\nmid u}}\frac {1}{||ua/q||}
\ll y+q
\end {eqnarray*}
the first sum in \eqref {split} is from part (A) and then \eqref {residue2}
\begin {eqnarray*}
%&&\sum _{u\leq y}d_k(u)\sum _{v\leq X/u}e\left (\frac {uva}{q}\right )=
\frac {X\mathcal N_{y}^{k-1}(q,q)}{q}+\mathcal O\left (\frac {Xq^{\theta -1}}{y^{\theta }}+y+q\right )
\end {eqnarray*}
whilst by assumption and then \eqref {wy2} the second is
\begin {eqnarray*}
%&&\sum _{v<X/y}\Big (D_{va/q}(X/v)-D_{va/q}(y)\Big )
%\\ \arrayfi =
%&&\sum _{v<X/y}\frac {1}{q_{v}}\sum _{r_1|\cdot \cdot \cdot |r_{k-2}|q_v\atop {\mathbf h|\mathbf r}}\frac {\mu (\mathbf h)}{\mathbf h}\left \{ \frac {X}{v}P_{k-1}(X/vq_v^2,\mathbf r)-yP_{k-1}(y,q_{v},\mathbf r)\right \} +\mathcal O\left (\frac {X}{y^\delta }\right )
%\\ \arrayfi =
%&&\frac {X}{q}\sum _{\mathbf h|\mathbf r|d|q}\frac {\mu (\mathbf h)}{\mathbf h}\sum _{H|d}\frac {\mu (H)P_{k+1}(X,y,q,\mathbf r,\mathbf h,d,H)}{H}
%\\ &&\hspace {30mm}+\hspace {4mm}\mathcal O\left (y+\sum _{v<X/y}\left \{ \frac {X}{v}\left (\frac {q}{X/v}\right )^{\Delta }+\frac {1}{\sqrt q}\left (\frac {q}{X/v}\right )^{\Delta /2}\right \} \right )
\frac {X}{q}\mathcal N_{X,y}^{k-1}(q,q)+\mathcal O\left (y+X\left (\frac {q}{y}\right )^\Delta +\frac {X}{\sqrt q}\left (\frac {q}{y}\right )^{\Delta /2}\right ).
\end {eqnarray*}
%so altogether we get a main term $X\mathcal N_{X,y}^{k-1}(q,q)/q$ and an error term
%\begin {eqnarray*}
%\ll \frac {Xq^{\theta -1}}{y^{\theta }}+\frac {X}{y}+y+q+X\left (\frac {q}{y}\right )^\Delta +\frac {X}{\sqrt q}\left (\frac {q}{y}\right )^{\Delta /2}.
%\end {eqnarray*}
Choosing $y=q\left (X/q\right )^{1/(1+\Delta )}$ gives a main term of the appropriate form and an error term
\[ \ll X\left (\frac {q}{X}\right )^{\Delta /(1+\Delta )}+\frac {X}{\sqrt q}\left (\frac {q}{X}\right )^{\Delta /2(1+\Delta )}\]
(it perhaps helping to note $\Delta /2\leq \theta $) and we have the claim for $k+1$.  \textbf {(C)} Write $\delta =(q,a)$ and $q=q'\delta $.  From part (B) the sum in question is
\begin {eqnarray*}
\frac {1}{q}\sum _{h|q}\sideset {}{'}\sum _{b=1}^he\left (-\frac {ab}{h}\right )\sum _{n\leq X}d_k(n)e\left (\frac {nb}{h}\right )=\frac {X}{q}\sum _{h|q}\frac {c_h(a)\mathcal M_X(h)}{h}+\mathcal O_q\left (X^{1-\Delta /2}\right )
\end {eqnarray*}
and here the sum in the main term is a polynomial in $\log X$ of degree $\leq k-1$ so comapring with \eqref {residue} the main term must be
\[ \frac {X}{\phi (q')}Res_{s=1}\Bigg \{ \frac {X^{s-1}}{s}\sum _{n=1\atop {(n,q)=\delta }}^\infty \frac {d_k(n)}{n^s}\Bigg \}%\hspace {10mm}\text{or}\hspace {10mm}\mathcal N_X^{k-2}(q,1)
\]
and therefore, from the sentence containing \eqref {gcd2}, we can assume $\delta =1$.  As $\mathcal N_X^{k-2}(q,1)/q$ is also a polynomial in $\log X$ of degree $\leq k-1$, it is enough to prove the claim with main term $X\mathcal N_X^{k-2}(q,1)/q$. For $k=2$ the claim holds by Theorem 1.1 of \cite {vaughandap} so we now suppose the claim is true for some $k\geq 2$ and show it also holds for $k+1$.  For a parameter $0<y\leq X$ we have 
\begin {eqnarray*}
\sum _{n\leq X\atop {n\equiv a(q)}}d_{k+1}(n)&=&\sum _{u\leq y\atop {(u,q)=1}}d_k(u)\sum _{v\leq X/u\atop {v\equiv a\overline {u}(q)}}1+\sum _{v<X/y\atop {(v,q)=1}}\sum _{y<u\leq X/v\atop {u\equiv a\overline {v}(q)}}d_k(u)
\\ &=&\frac {X}{q}\sum _{u\leq y\atop {(u,q)=1}}\frac {d_k(u)}{u}+\frac {1}{q}\sum _{v<X/y\atop {(v,q)=1}}\left \{ \frac {X\mathcal N_{X/v}^{k-2}(q,1)}{v}-y\mathcal N_y^{k-2}(q,1)\right \} 
%\\ &&\hspace {10mm}+\hspace {4mm}
+\mathcal O\left (y+\frac {q^{\Delta /2}X}{y^\Delta }+\frac {X}{y^{\delta}}\right ).
\end {eqnarray*}
From part (A), \eqref {residue2} and \eqref {wy3} we get a main term $X\mathcal N_{X,y}^{k-1}(q,1)/q$ and an error term
%\begin {eqnarray*}
%X\sum _{h_1,...,h_{k+1}|q}\frac {\mu (h_1)\cdot \cdot \cdot \mu (h_{k+1})}{h_1\cdot \cdot \cdot h_{k+1}}P_k(X,y,q,h_1,...,h_{k+1})
%\end {eqnarray*}
\begin {eqnarray*}
\ll \frac {X}{qy^{\theta  }}+y+\frac {q^{\Delta /2}X}{y^\Delta }+\frac {X}{y^{\delta}}.
\end {eqnarray*}
Note that 
\[ \theta \geq \delta \hspace {10mm}\frac {\delta }{1+\delta }=\frac {\Delta (1-\delta /2)}{1+\delta }=\left (1+\frac {\delta \Delta }{2(1+\delta)}\right )\frac {\delta }{1+\Delta }\]
so if $q\leq X^{\delta /(1+\delta )}$ choose $y=X^{1/(1+\delta )}$ and if $q\geq X^{\delta /(1+\delta  )}$ choose $y=\left (q^{\Delta /2}X\right )^{1/(1+\Delta )}$ - in both cases we get a total error 
\[ \ll \left (q^{\Delta /2}X\right )^{1/(1+\Delta )}+X^{1-\delta /(1+\delta )}\]
so we're done.  \textbf {(D)} From parts (A) and (C)%The main term in part (A) is a polynomial of degree $k$ in $\log X$ so from partial summation and comparing with part (C) we see that
\[ \sum _{n\leq X\atop {d|n}}d_k(n)=\frac {Xf_X(d,d)}{d}+\mathcal O\left (\left (\frac {X}{d}\right )^{1-\theta }\right )\]
so the sum in question is
\begin {eqnarray*}
%\sum _{d|q}\mu (q/d)d\sum _{n\leq X\atop {d|n}}d_k(n)&=&
X\sum _{D|d|q}\frac {\mu (q/d)\phi (D)\mathcal M_X(D)}{D}+\mathcal O\left (X^{1-\theta }q^{\theta }\right ).
\end {eqnarray*}
%and the sum here is
%\begin {eqnarray*}
%\sum _{D|q}\frac {\phi (D)\mathcal M_X(D)}{D}\sum _{d|q/D}\mu (q/dD)=\phi (q)\mathcal M_X(q).
%\end {eqnarray*}
\end {proof}
Next we use the previous lemma for our circle method input.
\begin {lemma}\label {barbandavenporthalberstam}
Let $\mathcal M_X(q)$ be as in Lemma \ref {kfold} and let
\[ \Delta _X(a/q)=\sum _{n\leq X}d_k(n)e(na/q)-\frac {X\mathcal M_X\left (q/(q,a)\right )}{q/(q,a)}.\]
Then for $1\leq Q,\gamma \leq x$ and errors $\mathcal O_\epsilon \left (x^\epsilon \cdot \cdot \hspace {1mm}\cdot \right )$
\begin {eqnarray*}
%V(x,Q)&\ll &xQ+\underbrace {x^{2-1/k}}_{k>2}
%\\ 
\sum _{q\leq \gamma }\frac {1}{q}\sideset {}{'}\sum _{a=1}^q\max _{X\leq x}|\Delta _X(a/q)|^2&\ll &x\gamma +\underbrace {x^{2-1/k}}_{k>2}
\\ \sum _{q\leq \gamma }\frac {1}{q^2}\sideset {}{'}\sum _{a=1}^q\max _{X\leq x}|\Delta _X(a/q)|^2&\ll &x+\underbrace {x^{2-8/(6k-3)}}_{k>2}.
\end {eqnarray*}
\end {lemma}
Letting $f_X(q,n)$ be as in Lemma \ref {kfold} we see that
\[ \frac {1}{q}\sum _{n=1}^qf_X(q,n)e\left (\frac {nb}{q}\right )=\frac {\mathcal M\left (q/(q,b)\right )}{q/(q,b)}\]
so letting $E_x(q,a)$ and $V(x,Q)$ be as in our theorems we have
%\[ E_x(q,a)=\frac {1}{q}\sum _{b=1}^qe\left (\frac {-ab}{q}\right )\left (\sum _{n\leq X}d_k(n)e(na/q)-\frac {\mathcal M\left (q/(q,b)\right )}{q/(q,b)}\right )\]
\begin {eqnarray}\label {clustog}
\sum _{a=1}^q|E_x(q,a)|^2&=&\frac {1}{q}\sum _{a=1}^q|\Delta _x(a/q)|^2
\end {eqnarray}
and therefore
\[ V(x,Q)\ll \sum _{q\leq Q}\frac {\log (Q/q)+1}{q}\sideset {}{'}\sum _{a=1}^q|\Delta _x(a/q)|^2\]
so the lemma implies Theorem \ref {largesieve}.
\begin {proof}
Let $X\leq x$ and drop $x^\epsilon $ and $\log x$ factors from error terms.  For $k=2$ Theorem 1.1 of \cite {blomer} says the LHS of \eqref {clustog} is $\ll X+q$ so let's assume $k>2$.  Take a parameter $0<R\leq \gamma $ and use \eqref {clustog} to write
\begin {alignat*}{4}
\sum _{q\leq \gamma }\frac {1}{q}\sideset {}{'}\sum _{a=1}^q|\Delta _X(a/q)|^2&\hspace {1mm}=\hspace {1mm}&&\left (\sum _{q\leq R}+\sum _{R<q\leq \gamma }\right )\frac {1}{q}\sideset {}{'}\sum _{a=1}^q|\Delta _X(a/q)|^2&\hspace {1mm}&=:A
\\ \sum _{q\leq \gamma }\frac {1}{q^2}\sideset {}{'}\sum _{a=1}^q|\Delta _X(a/q)|^2&\hspace {1mm}\leq \hspace {1mm}&&\sum _{q\leq R}\frac {1}{q}\sum _{a=1}^q|E_X(a/q)|^2+\sum _{R<q\leq \gamma }\frac {1}{q^2}\sideset {}{'}\sum _{a=1}^q|\Delta _X(a/q)|^2&\hspace {1mm}&=:B.\hspace {1mm}
\end {alignat*}
A large sieve inequality says
\[ \sum _{q\leq \gamma }\sideset {}{'}\sum _{a=1}^q|\Delta _X(a/q)|^2\ll X\gamma ^2+X^2\]
so from Lemma \ref {kfold} (B) 
\begin {eqnarray*}
A&\ll &X^{2-2\Delta }R^{1+2\Delta }+X^{2-\Delta }R^{\Delta }+X\gamma+\frac {X^2}{R}\ll X\gamma +X^{2-1/k} 
%\max _{R\leq Q}\Bigg \{ \frac {1}{R}\sum _{R\leq q\leq 2R}\sideset {}{'}\sum _{b=1}^q|\Delta _X(q,b)|^2\Bigg \} %\sum _{q\leq 2R}\frac {1}{q}\sideset {}{'}\sum _{b=1}^q|X_{q,b}|^2
%\\ &\ll &\max _{R\leq Q}\Bigg \{ \frac {1}{R}\sum _{q\sim R}\sum _{b=1}^q\left (|D_x(b/q)|^2+\left |\frac {x}{q}\sum _{N=1}^q\left (\sum _{d|q}c_d(b)f_x(d)\right )e\left (\frac {Nb}{q}\right )\right |^2\right )\Bigg \} .
%=:\max _{R\leq Q}\Upsilon _R
\end {eqnarray*}
on choosing $R=X^{1/k}$ whilst from Lemma \ref {kfold} (C) 
\begin {eqnarray*}
B&\ll &X^{2-2\Delta }R^{1+\Delta }+X^{2-2\delta }R+X+\frac {X^2}{R^2}\ll X^{2-4\delta /3}
%\max _{R\leq Q}\Bigg \{ \frac {1}{R}\sum _{R\leq q\leq 2R}\sideset {}{'}\sum _{b=1}^q|\Delta _X(q,b)|^2\Bigg \} %\sum _{q\leq 2R}\frac {1}{q}\sideset {}{'}\sum _{b=1}^q|X_{q,b}|^2
%\\ &\ll &\max _{R\leq Q}\Bigg \{ \frac {1}{R}\sum _{q\sim R}\sum _{b=1}^q\left (|D_x(b/q)|^2+\left |\frac {x}{q}\sum _{N=1}^q\left (\sum _{d|q}c_d(b)f_x(d)\right )e\left (\frac {Nb}{q}\right )\right |^2\right )\Bigg \} .
%=:\max _{R\leq Q}\Upsilon _R
\end {eqnarray*}  
on choosing $R=X^{2\delta /3}$ and perhaps noting $\delta \Delta /3\leq \Delta -\delta $.
\end {proof}
At one point we'll want to apply a variant of Perron's formula.
\begin {lemma}\label {k}
Take $\Phi :\mathbb N\rightarrow \mathbb C$ with
\[ \mathcal F(s)=\sum _{N=1}^\infty \frac {\Phi (N)}{N^s}\hspace {10mm}\text { absolutely convergent for }\sigma >1.\]
For $X\geq 1$ suppose that $K_X:\mathbb C\rightarrow \mathbb C$ satisfies:
\begin {eqnarray*}
\hspace {5mm}-&&K_X\text { is holomorphic except for possibly poles at $s\in \{ 0,-1,... \} $}
%\\ \hspace {5mm}-&&\text {there is a gap $\gg 1$ between the poles infinitely often}
%\\ \hspace {5mm}-&&\text {for each pole $Z$ of order $n$}\hspace {10mm}K_X(s)\ll \frac {\lambda _XX^{\sigma }}{|s|}\left (1+\frac {1}{|s-Z|^{n}}\right )\hspace {5mm}\text { for }s\not \in \{ 0,Z\} 
\\ \hspace {5mm}-&&\text {for all $R\in \mathbb N$ and $\sigma >-R-1$ with $s\not \in \{ 0,-1,...\} $}\hspace {10mm}K_X(s)\ll _\epsilon \frac {X^{\sigma +\epsilon }}{|s|}\sum _{r=1}^R\frac {X^{r}}{|(s+1)\cdot \cdot \cdot (s+r)|}
%\\ \hspace {5mm}-&&\text {away from the poles}\hspace {10mm}K(s)\ll \frac {\lambda X^{\sigma }}{|s|}
\\ \hspace {5mm}-&&\sum _{poles}Res \Bigg \{ \frac {K_X(s)}{N^s}\Bigg \} \hspace {5mm}\text { converges for any $N\in \mathbb N$}.
\end {eqnarray*}
Then for $X,Q\geq 1$ with $X/Q\not \in \mathbb Z$, $d\in \mathbb R$ and $c>1-d,0$
\begin {eqnarray*}
\sum _{N\leq X/Q}\frac {\Phi (N)}{N^d}\sum _{poles}Res \Bigg \{ \frac {K_X(s)}{(QN)^s}\Bigg \} =\int _{c\pm i\infty }\frac {\mathcal F(s+d)K_X(s)ds}{Q^s}.%+\mathcal O\left (\frac {\lambda E}{T^2}\right )
\end {eqnarray*}
%\begin {eqnarray*}
%E=(X/Q)^c\sum _{N=1}^\infty \frac {|\phi (N)|}{N^{c+d}}+(X/Q)^{1-d}\log (X/Q+2)\max _{X/2Q<N\leq 3X/2Q}|\Phi (N)|.
%\end {eqnarray*}
\end {lemma}
\begin {proof}
Take integer $R>X$ so that on the contour $\mathcal C$ consisting of the straight lines connecting
\[ c+iR\hspace {7mm}-R-1/2+iR\hspace {7mm}-R-1/2-iR\hspace {7mm}c-iR\hspace {7mm}\text { we have }\hspace {5mm}K_X(s)\ll \frac {X^{\sigma +\epsilon }e^X}{t}\]
%On these three lines we have in order the bounds
%\[ K(s)\ll \frac {\lambda X^{\sigma +1}}{T^2}\hspace {5mm}K(s)\ll \frac {\lambda }{Xt^2}\hspace {5mm}K(s)\ll \lambda X^{\sigma -2}\hspace {5mm}K(s)\ll \frac {\lambda X^{-R}e^{X}}{t}\]
and therefore
\begin {eqnarray*}
\int _\mathcal C\frac {K_X(s)ds}{(QN)^s}%\ll X\int _{c}^R\frac {(X/N)^\sigma ds}{T}\ll 
%\\ &&\hspace {10mm}\ll \hspace {4mm}
\ll X^\epsilon e^X\left (\frac {(X/QN)^{-R-1/2}+(X/QN)^c}{R|\log (X/QN)|}+(X/QN)^{-R}\log R\right )
%\ll \lambda e^X\left (\frac {(X/QN)^{c}+(X/QN)^{-R}}{T|\log (X/QN)|}+(X/QN)^{-R}\log (2T)\right )
\end {eqnarray*}
so by the Residue Theorem and pushing $R\rightarrow \infty $
\[ \int _{c\pm i\infty }\frac {K_X(s)}{(QN)^s}ds=\sum _{\text {poles}}Res \Bigg \{ \frac {K_X(s)}{(QN)^s}\Bigg \} \]
%\begin {eqnarray*}
%\int _{c\pm i\infty }\frac {K(s)}{(QN)^s}ds&=&\sum _{\text {poles}}Res \Bigg \{ \frac {K(s)}{(QN)^s}\Bigg \} +\mathcal O\left (\frac {\lambda (X/QN)^c}{T|\log (X/QN)|}\right )
%\end {eqnarray*}
as long as $QN<X$.  If $QN>X$ a similar argument but taking the rectangle to the right shows this last integral to be zero and we're done.
\end {proof}
%Let
%\[ W_{x,Q}(U)=F(t-UQ)dt\cdot \left \{ \begin {array}{ll}1&\text { if }U<x/Q\\ 0&\text { if }U>xQ\end {array}\right .\]
%and
%\[ w_{x,Q}(U)=W_{x,Q}(Ux/Q)\]
%so
%\[ w_{x,Q}\left (\frac {N}{x/Q}\right )=W_{x,Q}(N)=F(t-NQ)\cdot \left \{ \begin {array}{ll}1&\text { if }N<x/Q\\ 0&\text { if }N>xQ\end {array} \right .\]
%and
Now we apply the last result to what will be our case of interest. 
\begin {lemma}\label {w}
For $f:[1,\infty )\rightarrow \mathbb R$ integrable and $\Phi :\mathbb N\rightarrow \mathbb C$ with $\Phi (N)\ll _\epsilon N^\epsilon $ define for $X\geq 1$ and $\sigma >1$
%\[ X>0\hspace {6mm}s\in \mathbb C\hspace {6mm}\Phi :\mathbb N\rightarrow \mathbb C\]
\begin {eqnarray*}
%X\geq 1\hspace {2mm}\sigma >1\hspace {10mm}
F(X)&=&\int _1^Xf(t)dt\hspace {6mm}\text { and }\hspace {6mm}\mathcal F(s)=\sum _{N=1}^\infty \frac {\Phi (N)}{N^s}.
\end {eqnarray*}
If $f$ is smooth with $f^{(n)}(X)\ll _{n}X^{\epsilon -n}$ for all $n\geq 0$ then for $X-1,Q\geq 1$ and $c>0$
\begin {eqnarray*}
\sum _{u\leq (X-1)/Q}\frac {\Phi (u)}{u}F(X-uQ)&=&\int _{c\pm i\infty }\frac {\alpha (s+1)K_{X}(s)ds}{Q^s}%+\mathcal O\left (\frac {(t/Q)^c}{T}\right )
\end {eqnarray*}
where $K_X(s)$ satisfies the following properties:
\begin {eqnarray*}
\hspace {5mm}-&&K_X(s)\text { is holomorphic except for possibly poles at $s\in \{ 0,-1,... \} $}
%\\ \hspace {5mm}-&&Res_{s=-R}\{ K_X\} =\frac {(-1)^RF^{(R)}(X)}{R!}
%\\ K_X(s)&\ll &(X-1)^{s+1}\sum _{r=1}^{R+1}\frac {f^{(r-1)}(1)X^{r}}{|s|\cdot \cdot \cdot |s+r|}\hspace {6mm}\text {for }\sigma >-(R+1),\hspace {1mm}s\not \in -\mathbb N\hspace {5mm}R\in \mathbb N\cup \{ 0\} 
%\\ K_X(s)&\ll _\epsilon &\frac {(X-1)^{s+1}}{|s+1|}\text { for }s\not \in \{ |s+n|\leq \epsilon |n\geq 0\} 
%\\ \hspace {5mm}-&&\text {for each $n\geq 2$}\hspace {5mm}K_X(s)\ll \frac {(X-1)^{\sigma +1}}{|s(s+1)|}\left (1+\frac {1}{|s+n|}\right )\hspace {5mm}\text { for }s\not \in -\mathbb N\cup \{ 0\} 
\\ \hspace {5mm}-&&\text {for $R\in \mathbb N$ and $\sigma >-(R+1)$ with $s\not \in \{ 0,-1,...\} $}\hspace {10mm}K_X(s)\ll \sum _{r=1}^R\frac {(X-1)^{\sigma +r+\epsilon }}{|s\cdot \cdot \cdot (s+r)|}
%\frac {(X-1)^{\sigma +1}}{|s(s+1)|}
\\ \hspace {5mm}-&&\text {for $\sigma >-2$ with $s\not \in \{ 0,-1\} $}\hspace {5mm}\int _1^Xg(t)K_t(s)dt\ll \frac {X^{\sigma +2+\epsilon }}{|s(s+1)(s+2)|}\hspace {5mm}\text {for $g$ with }\left \{ \begin {array}{ll}g(X)\ll 1\\ g'(X)\ll X^{\epsilon -1}.\end {array}\right .
\end {eqnarray*}
Further, if $f(t)=(\log t)^D$ for $D\in \mathbb N\cup \{ 0\} $ then
\begin {eqnarray*}
\hspace {5mm}-&&\frac {K_X(s)}{Q^s}=\frac {X}{s}\sum _{n=0}^\infty P_n(X,Q)s^n\hspace {5mm}\text { for }|s|<1
\\ \hspace {5mm}-&&\frac {K_X(s)}{Q^s}=\frac {Q}{s+1}\sum _{n\geq 0}P_n(X,Q)(s+1)^n\hspace {5mm}\text { for }|s+1|<1
\end {eqnarray*}
where the $P_n(X,Q)$ are polynomials in $\log X,\log Q$ of degree $\leq D+n$.
%Test case: $f=1$.  Then $K_{t}(s)=t/s(s+1)$ and $F(t-Qu)=t-Qu$ and it says
%\[ \sum _{u\leq t/Q}f(u)(t-Qu)=t\int _{\pm c\pm iT}\frac {\alpha (s)(t/Q)^sds}{s(s+1)}\]
\end {lemma}
\begin {proof}
Define\footnote {should $K_X(s)$ come across as unmotivated then notice it is essentially the Mellin transform  
\[ \int _0^\infty w_X(z)z^{s-1}dz\hspace {5mm}\text { of }\hspace {5mm}w_X(z)=%\int _0^{X(1-z)}f(V)dV=
F(X(1-z))\left \{ \begin {array}{ll}1&\text { if }z<1-1/X\\ 0&\text { if }z>1+1/X\end {array}\right .\]}
for $\sigma >-1$ and $s\not =0$
\begin {eqnarray*}
%\alpha (s)=\sum _{N=1}^\infty \frac {\Phi (N)}{N^s}\hspace {6mm}\text { and }\hspace {6mm}
K_{X}(s)&=&%\int _0^\infty w_{x,Q}(U)U^{s-1}dU
%\\ &=&\int _0^\infty W_{x,Q}(Ux/Q)U^{s-1}dU
%\\ &=&\int _0^1F(t-xU)U^{s-1}dU
%\\ &=&\int _{0}^{t/x}F(t-xU)U^{s-1}dU
%\\ &=&\int _{0}^{t/x}\int _0^{t-xU}f(V)dVU^{s-1}dU
%\\ &=&\int _0^{t}f(V)\int _0^{(t-V)/x}U^{s-1}dUdV
%\\ &=&\frac {1}{sx^s}\int _0^{t}f(V)(t-V)^sdV
%\\ &=&
\frac {1}{s}\int _{0}^{X-1}f(X-W)W^sdW.
\end {eqnarray*}
Integrating any $R\geq 0$ times we get for $\sigma >-(R+1)$ and $s\not \in -\mathbb N\cup \{ 0\} $ 
\begin {eqnarray}\label {zoo}
K_X(s)&=&\sum _{r=1}^R\frac {f^{(r-1)}(1)(X-1)^{s+r}}{s\cdot \cdot \cdot (s+r)}
%\\ &&\hspace {0mm}+\hspace {4mm}
+\frac {1}{s\cdot \cdot \cdot (s+R)}\int _0^{X-1}f^{(R)}(X-W)W^{s+R}dW.
%\\ &\ll &(X-1)^{\sigma +1}\sum _{r=1}^{R+1}\frac {1}{|s|\cdot \cdot \cdot |s+r|}.\notag 
\end {eqnarray}
A classical argument using Morera's Theorem says that this integral is holomorphic for $\sigma >-(R+1)$ and the first two claimed properties of $K_x(s)$ follow, as does
\begin {equation}\label {resres}
Res_{s=-R}\{ K_X\} =\frac {F^{(R)}(X)}{(-R)!}.
\end {equation}
The third property follows since for $\sigma >-2$ with $s\not \in \{ 0,-1\} $
\begin {eqnarray*}
\int _1^Xg(t)K_t(s)dt&=&\frac {1}{s}\int _1^Xf(V)\int _{V}^Xg(t)(t-V)^sdtdV
\\ &=&\frac {1}{s(s+1)}\int _1^Xf(V)\left (g(X)(X-V)^{s+1}-\int _{V}^Xg'(t)(t-V)^{s+1}dtdV\right )
\\ &=&\frac {g(X)}{s(s+1)(s+2)}\left (f(1)(X-1)^{s+2}+\int _1^Xf'(V)(X-V)^{s+2}dV\right )-\frac {1}{s}\int _1^Xg'(t)K_t(s+1)dt
\end {eqnarray*}
which we can bound using the second property.  From \eqref {resres}
\begin {eqnarray*}
\sum _{\text {poles}}Res \Bigg \{ \frac {K_{X}(s)}{(Qu)^s}\Bigg \} &=&\int _0^{X-1}f(X-W)dW+\sum _{R=1}^\infty \frac {F^{(R)}(X)}{R!}\left (-Qu\right )^R=F(X-Qu)
\end {eqnarray*}
so the main claim follows from the first and second properties and from Lemma \ref {k}, and now we turn to the series expansions.  As
\[ X^{s+1}=\sum _{n\geq 0}\frac {(s+1)^n(\log X)^n}{n!}\hspace {6mm}\text { and }\hspace {6mm}-\frac {1}{s}=\sum _{n\geq 0}(s+1)^n\]
we get from \eqref {zoo} 
\begin {eqnarray*}
\frac {K_X(s)}{Q^{s}}%&=&\frac {Q}{s(s+1)}\left (f(1)\left (\frac {X-1}{Q}\right )^{s+1}+\int _0^{X-1}f'(X-W)\left (\frac {W}{Q}\right )^{s+1}dW\right )
%\\ 
&=&\frac {Q}{s(s+1)}\sum _{n\geq 0}\frac {(s+1)^n}{n!}\left (\underbrace {f(1)\left (\log \left (\frac {X-1}{Q}\right )\right )^n+\int _{0}^{X-1}f'(X-W)\left (\log \left (\frac {W}{Q}\right )\right )^ndW}_{\text {degree $\leq n+D$ polynomial}}\right )
%\\ &=&\frac {Q}{s(s+1)}\sum _{n\geq 0}\frac {(s+1)^n}{n!}P_n(\log X,\log Q)
\\ &=&\frac {Q}{s+1}\sum _{l\geq 0}(s+1)^l\sum _{n+m=l}\frac {P_n(\log X,\log Q)}{n!}
\end {eqnarray*}
and the expansion about $s=0$ is similar.
%\begin {eqnarray*}
%\frac {K_X(s)}{Q^s}&=&%\frac {1}{s}\int _{0}^{t}f(t-W)\left (W/Q\right )^sdW
%\\ &=&
%\frac {1}{s}\sum _{n\geq 0}\frac {s^n}{n!}\int _{0}^{X-1}f(X-W)\left (\log (W/Q)\right )^ndW.
%\\ &=&\frac {1}{s}\sum _{n\geq 0}\frac {s^n}{n!}P_n(X,Q)
%\end {eqnarray*}
\end {proof}
The remaining lemmas are concerned with writing the Dirichlet series arising from applying Perron's formula as ultimately a product of the Riemann zeta function and its derivatives.  The expression we want is Lemma \ref {arch2} and the inductive step for general $k$ is Lemma \ref {arch}. 
\begin {notation}\label {note}
Here we explain what we mean with the notation
\[ \mathcal O^*\hspace {5mm}d_k^w\hspace {5mm}f_k^{q}\hspace {5mm}\theta ^{}\hspace {5mm}
%\theta \hspace {5mm}
%\Delta ^{w,w'}\hspace {5mm}
\Delta \hspace {5mm}F_k^{s,w,w'}\hspace {5mm}F_k^{s}%\hspace {5mm}\lambda _s(w,w')%\hspace {5mm}\mathcal O^*\hspace {5mm}d_K^w
\]
which we use in Lemmas \ref {arch} and \ref {arch2}.  If $E_{s,w,w'}:\mathbb N\rightarrow \mathbb C$ 
\begin {eqnarray*}
&&-\hspace {2mm}\text {is holomorphic for $\sigma >0$ and differentiable for $w,w'\geq 0$}
\\ &&-\hspace {2mm}\text {is $>-1$ for $s$ real and positive}
\\ &&-\hspace {2mm}\text {is }\mathcal O(1/p^\sigma )\hspace {5mm}\text {for }\hspace {5mm}\text {primes $p$}\hspace {5mm}N\geq 1\hspace {5mm}\sigma >0\hspace {5mm}w,w'\geq 0
\end {eqnarray*}
then write $E_{s,w,w'}(p^N)=\mathcal O^*\left (1/p^\sigma \right )$.  For $k\in \mathbb N$ and $w\in \mathbb C$ define $d_k^w:\mathbb N\rightarrow \mathbb C$ through
\[ d_k^{w}(n)=\sum _{u_1\cdot \cdot \cdot u_k=n}\frac {1}{u_{k}^w}\]
and if there is multiplcative $f_k$ such that
\[ f_k\left (p^N\right )=d_k^w\left (p^N\right )+\mathcal O_{N,k}^*\left (1/p^\sigma \right )%\hspace {5mm}\text {for all primes $p$, all $N\geq 1$, and all $\sigma >0$}
\]
then say that $f_k$ is an $\mathcal S,\mathcal W$-function.  If $\Delta :\mathbb N\rightarrow \mathbb C$ is only squarefree kernel dependent and satisfies
\[ \Delta (p)=1+\mathcal O^*\left (1/p^\sigma \right ).%\hspace {5mm}\text {for all primes $p$ and all $\sigma >0$}.
\]
then say that $\Delta $ is an $\mathcal S,\mathcal W$-function.  If we repeat this discussion but without $w,w'$ parameters call $f_k$ and $\Delta $ just $\mathcal S$-functions.  Write $\theta $ for a function $\mathbb N\times \mathbb N\rightarrow \mathbb C$ which 
\begin {eqnarray*}
%&&-\hspace {2mm}\text {is holomorphic for $\sigma >0$ and differentiable for $w,w'\geq 0$}
%\\ 
&&-\hspace {2mm}\text {is positive and $\ll 1$ for real and positive $s$}
\\ &&-\hspace {2mm}\text {is multiplicative in both arguments}
%\\ &&-\hspace {2mm}\text {has $\theta (n,m)=1$ for $(n,m)=1$}
\\ &&-\hspace {2mm}\text {is only squarefree kernel dependent in the second argument}
\end {eqnarray*}
and then for given $\theta $ and multiplicative $f:\mathbb N\rightarrow \mathbb C$ write $f^q(n)=f(n)\theta (n,q)$.  
\\
\\ We make the convention that we use the same letters for all these functions even when the functions themselves may differ - their properties however remain intact.  For $N,N'\in \mathbb N$ write $\mathcal N=(N,N')$.  Write $F_k^{s,w,w'}(N,N')$ for a product of the form
\begin {eqnarray*}
%&&F_K^{s,\mathbf w}\left (\frac {\mathbf n}{(n,n')},(n,n')\right )=
%&&
f_{k}^{\mathcal N}\left (N/\mathcal N\right )f_{k}^{\mathcal N}\left (N'/\mathcal N\right )\Delta ^{}(\mathcal N)
\end {eqnarray*}
where each function is an $\mathcal S,\mathcal W$-function and similarly write $F_k^s(N,N')$ for the corresponding product with $\mathcal S$-functions. % Write $\lambda _s(w,w')$ for any quantity bounded and holomorphic for $\sigma >1/2$ and differentiable in $w,w'$ about zero.
\[ \hspace {110mm}\blacksquare \]
\end {notation}
\begin {lemma}\label {arch}
Assume the notation of Notation \ref {note} and define for $\sigma >1$ and $w,w'\geq 0$
\begin {eqnarray*}
\mathcal R_{N,N'}^{w,w'}(s)&=&\sum _{r,r'=1}^\infty \frac {F_k^s\left (rN,r'N'\right )}{[rN,r'N']^sr^wr'^{w'}}%f_{(rN,r'N')}^K\left (\frac {\mathbf r\mathbf N}{(rN,r'N')}\right )\Delta (rN,r'N')
\hspace {5mm}\text { and }\hspace {5mm}\mathcal H_{N,N'}^{w,w'}(s)=\sum _{h,h'=1}^\infty \frac {\mu (h)\mu (h')F_k^s\left ([h,N],[h',N']\right )}{hh'[[h,N],[h',N']]^sh^wh'^{w'}}.
\end {eqnarray*}
Then for $s=\sigma >1$ and $w,w'\geq 0$ 
\begin {eqnarray*}
&&\text {(A)}\hspace {10mm}\mathcal R_{N,N'}^{w,w'}(s)=\frac {\zeta (s+w+w')\zeta (s+w)^k\zeta (s+w')^k\mathcal \lambda _{w,w'}(s)F_{k+1}^{s,w,w'}\left (N,N'\right )}{[N,N']^s}
\\ &&\text {(B)}\hspace {10mm}\mathcal H_{N,N'}^{w,w'}(s)=\frac {\lambda _{w,w'}(s)F_{k}^{s,w,w'}\left (N,N'\right )}{[N,N']^s}
\end {eqnarray*}
where $\lambda _{w,w'}:\{ s\in \mathbb C|\sigma >1/2\} \rightarrow \mathbb C$ is holomorphic and $\ll 1$.
\end {lemma}
\begin {proof}
%For $d\in \mathbb N$ and a never-zero-multiplicative function $f$ define $f_d$ through
%\[ f(dn)=f(d)f(n)\prod _{p|n,d}\frac {f(p^{N+D})}{f(p^N)f(p^D)}=:f(d)f_d(n)\]
%where $p^N||n$, $p^D||d$, and for a never-zero-multiplicative-in-two-variables function $f$ define ${}_df(\cdot ,\cdot )$ through
%\[ f(dn,m)=f(d,m)f(n,m)\prod _{p|n,d}\frac {f(p^{N+D},m)}{f(p^N,m)f(p^D,m)}=:f(d,m){}_df(n,m)\]
%and similarly define $f_d(\cdot ,\cdot )$.  Note that $f_d(p^N,\cdot )=f_{p^D}(p^N,\cdot )$ and similarly for the others.  
%For $d\in \mathbb N$ and a never-zero-multiplicative function $f$ define $f_d$ to be the multiplicative function given by
%\text {for primes $p$ and $N\geq 0$}\hspace {5mm}
%f_d(p^N)=\frac {f(p^{N+D})}{f(p^D)}\hspace {10mm}\text {where }p^D||d\]
%which implies $f(dn)=f(d)f_d(n)$ for any $d,n\in \mathbb N$, implies $f_d(p^N)=f_1(p^N)$ if $p\nmid d$, and implies $f_d(p^N)=f_p(p^N)$ if $p|d$ and $f$ is only squarefree kernel dependent.  For a never-zero-multiplicative-in-two-variables function $f$ define ${}_df(\cdot ,\cdot )$ through ${}_df(n,m)=f(\cdot ,m)_d(n)$ and similarly define $f_d(\cdot ,\cdot )$.  %Note that if $f$ is a certain one of the three types of functions laid out in Notation \ref {note} then so is $f_d$ or ${}_df$ and we write
%\[ \{ f_K\} _d\hspace {5mm}{}_d\theta \hspace {5mm}\theta _d\hspace {5mm}\Delta _d\]
%for those functions.  Consequently for any $x,n\in \mathbb N$
Before turning to the claims themselves we need some preparatory work.  If $d\in \mathbb N$ appears in a $p$-factor of an Euler product then read $d$ as ``the highest power of $p$ which divides $d$".  For any $x,d\in \mathbb N$
\begin {eqnarray*}
\sum _{R=1}^\infty \frac {\theta (x,Rd)\Delta (Rd)}{R^{s+w+w'}}%&=&\theta (x,d)\Delta (d)\sum _{R=1}^\infty \frac {\theta _d(x,R)\Delta _d(R)}{R^{s+w+w'}}
%\\ &=&
=\prod _{p}\underbrace {\left (\sum _{R\geq 0}\frac {\theta (x,p^Rd)\Delta (p^Rd)}{p^{R(s+w+w')}}\right )}_{=:\mathfrak f_{x,d}}
\end {eqnarray*}
where $\mathfrak f_{1,1}=1+\mathcal O^*\left (1/p^{\sigma }\right )$.  This product is
\begin {eqnarray*}
%\prod _{p\nmid xd}\mathfrak f_{1,1}\prod _{p|xd}\mathfrak f_{x,d}=
\prod _p\mathfrak f_{1,1}\prod _{p|xd}\frac {\mathfrak f_{x,d}}{\mathfrak f_{1,1}}=:
\prod _{p}\mathfrak f_{1,1}\prod _{p|xd}\mathfrak g_{x,d}
\end {eqnarray*}
and here the infinite product is
\begin {eqnarray*}
%\prod _p\left (\sum _{R\geq 0}\frac {\Delta (p^R)}{p^{R(s+w+w')}}\right )&=&
\zeta (s+w+w')\prod _{p}\left (1+\mathcal O^*\left (\frac {1}{p^{2\sigma +w+w'}}\right )\right )
%\\ &=&
%=\zeta (s+w+w')\lambda _s(w,w')
\end {eqnarray*}
and the finite product
\begin {eqnarray*}
%\prod _{p|x\atop {p\nmid d}}\mathfrak g_{x,1}\prod _{p|d\atop {p\nmid x}}\mathfrak g_{1,p}\prod _{p|x,d}\mathfrak g_{x,p}=
\prod _{p|x}\mathfrak g_{x,1}\prod _{p|d}\mathfrak g_{1,p}\prod _{p|x,d}\frac {\mathfrak g_{x,p}}{\mathfrak g_{x,1}\mathfrak g_{1,p}}=f_1^{}(x)\Delta ^{}(d)\theta ^{}(x,d)
\end {eqnarray*}
so for any $x,d\in \mathbb N$
\begin {eqnarray}\label {helpu}
&&\sum _{R=1}^\infty \frac {f_k^{Rd}(x)\Delta (Rd)}{R^{s+w+w'}}=\zeta (s+w+w')\lambda _{w,w'}(s)f_k^{}(x)\theta ^{}(x,d)\Delta ^{}(d).\hspace {10mm}
\end {eqnarray}
We also have
\begin {eqnarray*}
\sum _{H,h,h'=1\atop {(h,h')=1\atop {(hh',q)=1\atop {(H,hh'q)=1}}}}^\infty \frac {\Delta (H)f_k(h)f_k(h')}{H^{2+s+w+w'}h^{1+s+w}h'^{1+s+w'}}&=&\lambda _{w,w'}(s)\sum _{h,h'=1\atop {(h,h')=1\atop {(hh',q)=1}}}^\infty \frac {\Delta (hh'q)f_k(h)f_k(h')}{h^{1+s+w}h'^{1+s+w'}}
\\ &=&\lambda _{w,w'}(s)\Delta ^{}(q)\sum _{n=1\atop {(n,q)=1}}^\infty \frac {f_k^{}(n)}{n^{1+s+w+w'}}%\sum _{hh'=n\atop {(h,h')=1}}^\infty \frac {1}{\mathbf h^{w}}s
\end {eqnarray*}
so that a similar calculation to that just done leads to
\begin {eqnarray}\label {helpu2}
&&\sum _{H,h,h'=1\atop {(h,h')=1\atop {(hh',q)=1\atop {(H,hh'q)=1}}}}^\infty \frac {\Delta (H)f_k(h)f_k(h')}{H^{2+s+w+w'}h^{1+s+w}h'^{1+s+w'}}=\lambda _{w,w'}(s)\Delta ^{}(q).
\end {eqnarray}
Now we turn to the claims of the lemma.  Any boldface letter, say $\mathbf v$, will be understood to mean a two-dimensional vector, written $\mathbf v=(v,v')$, and we then write $\mathbf v'=(v',v)$.  %For any set $X$ we write $\mathbf v\in X$ to mean $v,v'\in X$.  
A function with a vector appearing in its argument, say $f(\mathbf v)$, is to be understood to mean $f(v)f(v')$.  A condition $(\mathbf a,\mathbf b)=1$ will mean $(a,b)=(a',b')=1$ and $\mathbf d|\mathbf n$ will mean $d|n$ and $d'|n'$.  Write $\mathcal N=(N,N')$.  \textbf {(A)} For any $a,b\in \mathbb N$ 
\[ (ab,a'b')=(a,a')\left (\frac {a}{(a,a')},\frac {b'}{(b,b')}\right )\left (\frac {b}{(b,b')},\frac {a'}{(a,a')},\right )(b,b')\]
so using \eqref {helpu} and writing $\mathbf M=\mathbf N/\mathbf n\mathcal N$ %Define(
%\begin {eqnarray*}
%\mathcal W_s^{\mathbf q}(w,w')&=&\frac {1}{N^wN'^{w'}}\sum _{r,r'=1\atop {(\mathbf r,\mathbf q)=1}}^\infty \frac {1}{[rN,r'N']^sr^wr'^{w'}}%f_{(rN,r'N')}^K\left (\frac {\mathbf r\mathbf N}{(rN,r'N')}\right )\Delta (rN,r'N')
%F_K^\mathbf w\left (\frac {\mathbf r\mathbf N}{(rN,r'N')},(rN,r'N')\right ).
%\end {eqnarray*}
%We have
%\begin {eqnarray}\label {bb}
%\mathcal W_s^{\mathbf q}(w,z)&=&\frac {1}{[N,N']^s}\sum _{\mathbf n|\mathbf N/\mathcal N\atop {(\mathbf n,\mathbf q')=1}}\frac {1}{n^wn'^z}\sum _{R,r,r'=1\atop {(r,r')=1\atop {(\mathbf r,\mathbf q\mathbf n\mathbf M)=1\atop {(R,qq')=1}}}}^\infty \frac {1}{R^{s+w+z}(rr')^sr^wr'^z}f\left (\mathbf r\mathbf M\right )\theta \left (\mathbf r\mathbf M,Rnn'\mathcal N\right )\Delta (Rnn'\mathcal N)\notag 
%\\ &=&\frac {\zeta (s+w+z)\mathcal F(s+w+z)\Delta (qq')}{[N,N']^s}\notag 
%\\ &&\hspace {0mm}\sum _{\mathbf n|\mathbf N/\mathcal N\atop {(\mathbf n,\mathbf q')=1}}\frac {\Delta (nn'\mathcal N)\Delta (nn'\mathcal N,qq')}{n^zn'^w}\sum _{r,r'=1\atop {(r,r')=1\atop {(\mathbf r,\mathbf q\mathbf n\mathbf M')=1}}}^\infty \frac {1}{(rr')^sr^wr'^z}f\left (\mathbf r\mathbf M\right )\theta \left (\mathbf r\mathbf M,nn'\mathcal N\right )
%\end {eqnarray}
\begin {eqnarray*}
\mathcal R_{N,N'}^{w,w'}(s)&=&\frac {1}{[N,N']^s}\sum _{\mathbf n|\mathbf N/\mathcal N}\frac {1}{{n'}^{w}n^{w'}}\sum _{R,r,r'=1\atop {(r,r')=1\atop {(\mathbf r,\mathbf n\mathbf M')=1}}}^\infty \frac {f_{k}^{Rnn'\mathcal N}\left (\mathbf r\mathbf M\right )%\theta \left (\mathbf r\mathbf M,Rnn'\mathcal N\right )
\Delta (Rnn'\mathcal N)}{R^{s+w+w'}(rr')^sr^wr'^{w'}}
\\ &=&\frac {\zeta (s+w+w')\lambda _s(w,w')}{[N,N']^s}\sum _{\mathbf n|\mathbf N/\mathcal N}\frac {\Delta ^{}(nn'\mathcal N)}{{n'}^wn^{w'}}\sum _{r,r'=1\atop {(r,r')=1\atop {(\mathbf r,\mathbf n\mathbf M')=1}}}^\infty \frac {f_{k}^{nn'\mathcal N}\left (\mathbf r\mathbf M\right )%\theta ^{}\left (\mathbf r\mathbf M,nn'\mathcal N\right )
}{(rr')^sr^wr'^{w'}}
\\ &=&\frac {\zeta (s+w+w')\lambda _s(w,w')}{[N,N']^s}\Delta ^{}(\mathcal N)
%\\ &&\hspace {0mm}
\sum _{\mathbf n|\mathbf N/\mathcal N}\frac {\Delta ^{}_\mathcal N(\mathbf n)}{n'^wn^{w'}}\sum _{l=1}^\infty \frac {a_{\mathbf n,\mathbf M,nn'\mathcal N}(l)}{l^s}
\end {eqnarray*}
where
\[ a_{\mathbf n,\mathbf M,nn'\mathcal N}^{}(l)=\sum _{rr'=l\atop {(r,r')=1\atop {(\mathbf r,\mathbf n\mathbf M')=1}}}\frac {f_k^{nn'\mathcal N}(\mathbf r\mathbf M)}{r^wr'^{w'}}
.\]
Make the same convention for $M,M'$ as for $d$ in the second sentence of this proof.  The $l$ series above is
\begin {eqnarray*}
%\sum _{l=1}^\infty \frac {a_{\mathbf n,\mathbf M,nn'\mathcal N}(l)}{l^s}&=&
\prod _{p}\underbrace {\left (\sum _{l\geq 0}\frac {a_{\mathbf n,\mathbf M,nn'\mathcal N}(p^l)}{p^{ls}}\right )}_{=:\mathfrak a_{\mathbf n,\mathbf M,nn'\mathcal N}}
%\\ &=:&
%=:\prod _p\mathfrak a_{\mathbf n,\mathbf M,nn'\mathcal N}^{}
%\\ &=&
=\prod _{p|NN'}\underbrace {\frac {\mathfrak a_{\mathbf n,\mathbf M,nn'\mathcal N}}{\mathfrak a_{\mathbf 1}}}_{=:\mathfrak b_{\mathbf n,\mathbf M,nn'\mathcal N}}\underbrace {\prod _{p}\mathfrak a_{\mathbf 1}}_{=:\mathcal A(w,w')}
\end {eqnarray*}
and here the finite product is %(split the $p$ first according to whether $p|NN'/\mathcal N^2$ or not)
\begin {eqnarray*}
%&&\prod _{p|N/\mathcal N}\mathfrak b_{n,1,N/n\mathcal N,1,n\mathcal N}\prod _{p|N'/\mathcal N}\mathfrak b_{1,n',1,N'/n'\mathcal N,n'\mathcal N}\prod _{p|NN'\atop {p\nmid NN'/\mathcal N^2}}\mathfrak b_{1,1,1,1,\mathcal N}
%\\ \arrayfi =
\prod _{p|N/\mathcal N}\mathfrak b_{n,1,N/n\mathcal N,1,n\mathcal N}\prod _{p|N'/\mathcal N}\mathfrak b_{1,n',1,N'/n'\mathcal N,n'\mathcal N}\prod _{p|\mathcal N\atop {p\nmid NN'/\mathcal N^2}}\mathfrak b_{1,1,1,1,p}%\prod _{p|\mathcal N,NN'/\mathcal N^2}\frac {1}{\mathfrak b_{1,1,1,1,p}}
%\\ &&\hspace {15mm}%\prod _{p|\mathcal N,N'/\mathcal N}\frac {1}{\mathfrak b_{1,1,1,1,p}}
%\\ \arrayfi =
%\\ &&\hspace {25mm}=:\hspace {4mm}
=:b_{\mathbf n,\mathcal N}\left (\mathbf N/\mathcal N\right )b_{NN'/\mathcal N^2}(\mathcal N)
%\\ \arrayfi =f_\mathbf n(\mathbf N)\theta _{\mathbf n}(\mathbf N,\mathcal N)\Delta _{\mathbf n,\mathbf M}(\mathcal N)
\end {eqnarray*}
so
\begin {eqnarray}\label {2ogloch}
\mathcal R_{N,N'}^{w,w'}(s)&=&\frac {\zeta (s+w+w')\lambda _{w,w'}(s)\mathcal A(w,w')}{[N,N']^s}b_{NN'/\mathcal N^2}(\mathcal N)B_{\mathcal N}^{}\left (\mathbf N/\mathcal N\right )\hspace {10mm}
%\\ &&\hspace {0mm}\sum _{\mathbf n|\mathbf N/\mathcal N}\frac {\Delta (\mathbf n)f(\mathbf M)\theta (\mathbf M,\mathbf n\mathcal N)\Delta (\mathbf n,\mathcal N)}{n^wn'^z}f(\mathbf N/\mathcal N,\mathbf n).
\end {eqnarray}
where
\begin {eqnarray}\label {hwna}
B_q^{}(N)&:=&\sum _{n|N}\frac {\Delta ^{}_q(n)f_k^{nq}(N/n)%\theta ^{}(N/n,nq)
b_{n,q}(N)}{n^{w'}}\notag 
%\\ &=:&
=:\sum _{n|N}\underbrace {\lambda _{N,q}(n)b_{n,q}(N)}_{\text {multiplicative}}\notag 
%\\ &=&
=B_1^{}(N)\prod _{p|N,q}\frac {B_p^{}(p^N)}{B_1^{}(p^N)}\notag 
%\\ &=&B_1^{}(N)\prod _{p|N,q}\left (1-\frac {B_p(p^N)-B_1(p^N)}{B_1^{}(p^N)}\right ).
\end {eqnarray}
%Note that for any $d,l\in \mathbb N$
%\[ \{ f_K^s\} _d(p^l)\ll l\hspace {6mm}\]
%For $\sigma >1/2$ 
%\begin {equation}\label {pp}
%\sum _{l\geq 1}\frac {d_k(p^l)}{p^{ls}}=\sum _{l\geq 1}{k+l-1}\choose {k-1}\frac {1}{p^{ls}}\ll \frac {1}{p^\sigma }
%\end {equation}
so we're done if we show
\begin {eqnarray}\label {aaa}
\mathcal A_s(w,w')=\zeta (s+w)^k\zeta (s+w')^k\lambda _{w,w'}(s)\hspace {7mm}\mathfrak b_{1,1,1,1,p}=1+\mathcal O^*\left (p^{-\sigma }\right )\hspace {7mm}B_q(N)=f_k^q(N).
\end {eqnarray}
For $N$ a power of $p$ and any $q\in \mathbb N$
%\begin {eqnarray*}
%\frac {f_{K}^{nq}\left (Np^l/n\right )%\theta ^{}\left (Np^{l}/n,nq\right )
%}{n^{w'}}&=&\frac {d_K\left (Np^l/n\right )}{n^{w'}}\left (1+\mathcal O\left (\frac {1}{p^\sigma }\right )\right )
%\end {eqnarray*}
\begin {eqnarray*}
\sum _{n|N}\lambda _{N,q}(n)&=&d_{k+1}^{w'}(N)+O^*\left (p^{-\sigma }\right )
\end {eqnarray*}
and for $l\geq 1$ %and $n|p^N$
\begin {eqnarray*}
%\theta ^{}(N/n,nq)
a_{n,1,N/n,1,nq}(p^l)&=&\sum _{rr'=p^l\atop {(r,r')=1\atop {(r,n)=1\atop {(r',N/n)=1}}}}\frac {f_{k}^{nq}\left (rN/n\right )f^{nq}(r')%\theta ^{}\left (Nr/n,nq\right )\theta ^{}\left (r',nq\right )
}{r^wr'^{w'}}
%\\ &=&
=\underbrace {\frac {f_K^{nq}\left (Np^l/n\right )%\theta ^{}\left (Np^{l}/n,nq\right )\theta ^{}(1,nq)
}{p^{lw}}}_{p\nmid n}+\underbrace {\frac {f_{k}^{}\left (N/n\right )}{p^{lw'}}}_{p\nmid N/n}
%\\ \arrayfi \ll d_K(Np^l)%+d_K(p^l)
%\\ \arrayfi =d_K^s(Np^l)
%\\ \arrayfi =\underbrace {\frac {d_K(Np^{l})}{p^{lw}}+\mathcal O\left (\frac {d_K(Np^l)}{p^s}\right )}_{p\nmid n}+\underbrace {\mathcal O(l)}_{p\nmid N/n}
%\\ \arrayfi \ll \left \{ \begin {array}{ll}&\text { if }\\ &\text { if }\end {array}\right .
\end {eqnarray*}
so
\begin {eqnarray*}
\sum _{l\geq 1}\frac {1}{p^{ls}}\sum _{n|N}\lambda _{N,q}(n)a_{n,1,N/n,1,nq}(p^l)&=&%\sum _{n|N}f_K^{nq}(N/n)
%\theta (N/n,nq)
\sum _{l\geq 1}\frac {1}{p^{ls}}\left (\frac {\Delta _q(1)f_k^{q}\left (Np^l\right )%\theta ^{}\left (Np^{l},q\right )
}{p^{lw}}+\frac {\Delta _q(N)f^{Nq}(p^l)%\theta ^{}\left (p^l,Nq\right )
}{p^{lw'}}\right )
%\\ &=&
=\mathcal O_{N,k}^*\left (p^{-\sigma }\right )
\end {eqnarray*}
%\begin {eqnarray*}
%&&\lambda _{N,q}(n)a_{n,1,N/n,1,nq}(p^l)
%\\ \arrayfi =\Delta ^{w,w'}(n)\Delta ^{w,w'}(n,q)\left (\underbrace {\frac {1}{p^{lw}}f_K^{w,w'}\left (Np^l/n\right )\theta ^{w,w'}\left (Np^{l}/n,nq\right )}_{p\nmid n}+\underbrace {\frac {1}{p^{lw'}}f_K^{w,w'}\left (p^l\right )\theta ^{w,w'}\left (p^l,nq\right )}_{p\nmid N/n}\right )
%\\ \arrayfi \ll d_K(Np^l)+d_K(p^l)
%\\ \arrayfi =d_K^s(Np^l)
%\\ \arrayfi =\underbrace {\frac {d_K(Np^{l})}{p^{lw}}+\mathcal O\left (\frac {d_K(Np^l)}{p^s}\right )}_{p\nmid n}+\underbrace {\mathcal O(l)}_{p\nmid N/n}
%\\ \arrayfi \ll \left \{ \begin {array}{ll}&\text { if }\\ &\text { if }\end {array}\right .
%\end {eqnarray*}
%so
%\begin {eqnarray*}
%\sum _{l\geq 1}\frac {1}{p^{ls}}\sum _{p|n|p^N}f^s_K(N/n)\theta ^s(N/n,nq)a_{n,1,p^N/n,1,nq}(p^l)\ll \frac {1}{p^s}
%\end {eqnarray*}
and we conclude that for $N$ a power of $p$ and for any $q\in \mathbb N$ 
\begin {eqnarray}\label {wcw}
\mathfrak a_\mathbf 1B_q(N)&=&%\sum _{n|N}\lambda _{N,q}(n)+\sum _{l\geq 1}\frac {1}{p^{ls}}\sum _{n|N}\lambda _{N,q}(n)a_{n,1,N/n,1,nq}(p^l)
%\\ &=&
d_{k+1}^{w'}(N)+O_{N,k}^*\left (p^{-\sigma }\right ).
\end {eqnarray}
For $l$ a power of $p$
\begin {eqnarray*}
a_{1,1,1,1,q}(l)&=&\left (\frac {f_k^q(l)}{l^w}+\frac {f^q(l)}{l^{w'}}\right )
%\\ &=&=d_K(l)\left (\frac {1}{l^w}+\frac {1}{l^{w'}}\right )+\mathcal O\left (\frac {d_K^{w'}(l)}{p^\sigma }\right )
\hspace {10mm}\text {so}\hspace {10mm}
\sum _{l\geq 2}\frac {a_{1,1,1,1,q}(p^l)}{p^{ls}}=\mathcal O_k^*\left (\frac {1}{p^{2\sigma }}\right )%\sum _{l\geq 2}\frac {1+p^{-l(w+w')}}{p^{ls}}\ll \frac {1}{p^{2s}}\frac {1}{1-p^{-s}}+\frac {1}{p^{2(s+w+w')}}\frac {1}{1-p^{-(s+w+w')}}\ll \frac {1}{p^{1+\delta }}\]
%and
%\[ a_\mathbf 1(p)=f_1(p)\left (\frac {1}{p^w}+\frac {1}{p^{w'}}\right )\]
\end {eqnarray*}
so
\begin {eqnarray*}
\mathfrak a_{1,1,1,1,q}&=&1+\frac {k}{p^s}\left (\frac {1}{p^w}+\frac {1}{p^{w'}}\right )+\mathcal O_k^*\left (\frac {1}{p^{2\sigma }}\right )\notag 
\\ &=&\left (1+\frac {1}{p^{s+w}}\right )^k\left (1+\frac {k}{p^{s+w'}}\right )^k+\mathcal O_k^*\left (\frac {1}{p^{2\sigma }}\right )
\\ &=&1+\mathcal O_k^*\left (\frac {1}{p^\sigma }\right ).\notag 
\end {eqnarray*}
The second equality here implies the first equality of \eqref {aaa}, the third equality here implies the second equality of \eqref {aaa}, and the third equality here with \eqref {wcw} implies the third equality of \eqref {aaa}.  \textbf {(B)} The sum in question is
\begin {eqnarray*}
&&\left (\sum _{\mathbf j|\mathbf N}\frac {\mu (\mathbf j)}{\mathbf j^\mathbf w}\right )\sum _{h,h'=1\atop {(\mathbf h,\mathbf N)=1}}^\infty \frac {\mu (\mathbf h)F(hN,h'N')}{\mathbf h[hN,h'N']^sh^wh'^{w'}}
%\\ \arrayfi =\left (\sum _{\mathbf j|\mathbf N}\frac {\mu (\mathbf j)}{j'^{w+1}j^{w'+1}}\right )\sum _{\mathbf h=1\atop {(\mathbf h,\mathbf N)=1}}^\infty \frac {\mu (\mathbf h)f_K(\mathbf h\mathbf M)\theta (\mathbf h\mathbf M,Hnn'\mathcal N)\Delta (Hnn'\mathcal N)}{\mathbf h[hN,h'N']^sh^wh'^{w'}}
\end {eqnarray*}
and here the $\mathbf h$ sum is, as in (A) and writing $\mathbf M=\mathbf N/\mathbf n\mathcal N$,
\begin {eqnarray*}
&&\frac {1}{[N,N']^s}\sum _{\mathbf n|\mathbf N/\mathcal N\atop {(\mathbf n,\mathcal N)=1}}\frac {1}{\mathbf n^{1+\mathbf w'}}\sum _{H,h,h'=1\atop {(h,h')=1\atop {(\mathbf h,\mathbf M'\mathbf N)=1\atop {(H,\mathbf N)=1}}}}^\infty \frac {\mu (\mathbf h\mathbf n'H)f_k(\mathbf h\mathbf M)\theta (\mathbf h\mathbf M,Hnn'\mathcal N)\Delta (Hnn'\mathcal N)}{H^{2+s+w+w'}\mathbf h^{1+s+\mathbf w}}
\\ \arrayfi =\frac {\Delta (\mathcal N)}{[N,N']^s}\left (\sum _{\mathbf n|\mathbf N/\mathcal N\atop {(\mathbf n,\mathcal N)=1}}\frac {\mu (\mathbf n)f_k(\mathbf M)\theta (\mathbf M,\mathbf n\mathcal N)\Delta (\mathbf n)}{\mathbf n^{1+\mathbf w'}}\right )\sum _{H,h,h'=1\atop {(h,h')=1\atop {(hh',NN')=1\atop {(H,hh'NN')=1}}}}^\infty \frac {\mu (\mathbf h)\mu (H)^2f_k(\mathbf h)\Delta (H)}{H^{2+s+w+w'}\mathbf h^{1+s+\mathbf w}}.
\end {eqnarray*}
The $j$ sum is $\Delta (\mathbf N)$, the $H,\mathbf h$ sum is $\lambda _{w,w'}(s)\Delta (NN')$ from \eqref {helpu2}, and the $\mathbf n$ sum is $F_\mathcal N(\mathbf N/\mathcal N)$ where
\begin {eqnarray*}
F_q(N)&=&\sum _{n|N\atop {(n,q)=1}}\frac {\mu (n)f_k(N/n)\theta (N/n,nq)\Delta (n)}{n^{1+w'}}=F_1(N)\prod _{p|N,q}\frac {F_p(p^N)}{F_1(p^N)}
\end {eqnarray*}
so putting everything together the sum in question is
\begin {eqnarray*}
\frac {\lambda _{w,w'}(s)\Delta (\mathbf N)\Delta (NN')\Delta (\mathcal N)F_\mathcal N(\mathbf N/\mathcal N)}{[N,N']^s}
\end {eqnarray*}
and we're done on noting
\[ F_q(p^N)=d_k^{w'}(p^N)+\mathcal O_{N,k}^*\left (p^{-1}\right ).\]
\end {proof}
\begin {lemma}\label {arch2}
Let $P(\cdot ,...,\cdot )$ be a polynomial of degree $D$ in the variables $\log (\cdot ),...,\log (\cdot )$, take $k\in \mathbb N$, and write $Z=k(k+2)$.  Then for $\sigma >1$
\begin {eqnarray*}
\sum _{q=1}^\infty \frac {\phi (q)/q}{q^s}\sum _{r_k|\cdot \cdot \cdot |r_1|q\atop {r_k'|\cdot \cdot \cdot |r_1'|q\atop {h_1|r_1,...,h_k|r_k\atop {h_1'|r_1',...,h_k'|r_k'}}}}\frac {\mu (h_1)\mu (h_1')\cdot \cdot \cdot \mu (h_k)\mu (h_k')}{h_1h_1'\cdot \cdot \cdot h_kh_k'}P(q,\mathbf r,\mathbf r',\mathbf h,\mathbf h')
%\\ &&\hspace {10mm}=\hspace {4mm}
=\sum _{k_i\geq 0\atop {k_0+\cdot \cdot \cdot +k_Z\leq D}}\lambda _{\mathbf k}(s)\zeta ^{(k_1)}(s)\cdot \cdot \cdot \zeta ^{(k_Z)}(s)
\end {eqnarray*}
where for $\sigma >1/2$ the $\lambda _\mathbf k(s)$ are holomorphic and $\ll 1$.  The claim remains true for $k=0$ if we remove the $r_1,...,h_k'$ sum and variables.
\end {lemma}
\begin {proof}
%Before we start the proof let us work out cleanly derivatives involving powers of $\zeta (s)$.  
Write $\lambda (s)$ for a quantity holomorphic and $\ll 1$ for $\sigma >1/2$, let $\Delta ,F_K^s$ be as given in Notation \ref {note}, and assume w.l.o.g. that $s=\sigma $.  From \eqref {helpu}
\begin {eqnarray*}
\sum _{q=1}^\infty \frac {\Delta (qR)}{q^s}&=&\lambda (s)\zeta (s)\Delta (R)
\end {eqnarray*}
so the sum in question is
\begin {eqnarray}\label {spark}
%&&\sum _{0\leq d\leq D}\lambda _d(s)\zeta ^{(D-d)}(s)\sum _{r_1|\cdot \cdot \cdot |r_k\atop {r_1'|\cdot \cdot \cdot |r_k'\atop {h_1|r_1,...,h_k|r_k\atop {h_1'|r_1',...,h_k'|r_k'}}}}\frac {\mu (h_1)\cdot \cdot \cdot \mu (h_k')\Delta ([r_k,r_k'])}{h_1\cdot \cdot \cdot h_k'[r,r']^s}P_d(\mathbf r,\mathbf r',\mathbf h,\mathbf h')\notag 
%\\ \arrayfi =
\sum _{0\leq d\leq D}\lambda _d(s)\zeta ^{(D-d)}(s)\sum _{r_k|\cdot \cdot \cdot |r_1\atop {r_k'|\cdot \cdot \cdot |r_1'\atop {h_1|r_1,...,h_k|r_k\atop {h_1'|r_1',...,h_k'|r_k'}}}}\frac {\mu (h_1)\cdot \cdot \cdot \mu (h_k')F_1^s(r_k,r_k')}{h_1\cdot \cdot \cdot h_k'[r,r']^s}P_d(\mathbf r,\mathbf r',\mathbf h,\mathbf h')\hspace {10mm}
\end {eqnarray}
for some polynomials $P_d(\cdot ,...,\cdot )$ in $\log (\cdot ),...,\log (\cdot )$ of degree $d$.  For $n,d\in \mathbb N\cup \{ 0\} $ write $P^{n,N}(s)$ for a quantity of the form
\[ \sum _{k_i\geq 0\atop {k_1+\cdot \cdot \cdot +k_n\leq N}}\lambda _\mathbf k(s)\zeta (s)^{(k_1)}\cdot \cdot \cdot %\zeta ^{(k_K)}(s)\zeta (s)^{(k_1')}\cdot \cdot \cdot 
\zeta ^{(k_n)}(s)\]
so that
\begin {eqnarray}\label {derivative}
&&\frac {d}{dw'^{R'}w^R}\Bigr |_{\mathbf w=\mathbf 0}\lambda _{w,w'}(s)\zeta (s+w+w')\zeta (s+w)^K\zeta (s+w')^K
%\\ &&\hspace {35mm}
=P^{2K+1,R+R'}(s).\hspace {10mm}
\end {eqnarray}
Define $\mathcal R_{N,N'}^{w,w'}(s)$ as in Lemma \ref {arch}, define
\[ \mathcal R_{N,N'}(s)=\sum _{r,r'=1}^\infty \frac {(\log r)^R(\log r')^{R'}F_K\left (rN,r'N'\right )}{[rN,r'N']^s}\hspace {5mm}\text { so that }\hspace {6mm}\frac {d}{dw^Rdw'^{R'}}R_{N,N'}^{w,w'}(s)\Bigr |_{\mathbf w=\mathbf 0}=\mathcal R_{N,N'}(s)
\]
and define
\begin {eqnarray}\label {blodyn0}
\mathcal F_{\mathbf N}(s)%&:=&\sum _{r,r',h,h'=1\atop {[h,N]|r\atop {[h',N']|r'}}}^\infty \frac {\mu (h)\mu (h')(\log r)^R(\log r')^{R'}(\log h)^H(\log h')^{H'}}{hh'[r,r']^s}F_K\left (\frac {\mathbf r}{(r,r')},(r,r')\right )
%\\ 
&=&\sum _{h,h'=1}^\infty \frac {\mu (h)\mu (h')(\log h)^H(\log h')^{H'}\mathcal R_{[h,N],[h',N']}(s)}{hh'}.
%\\ &=&\sum _{r,r',h,h'=1}^\infty \frac {\mu (h)\mu (h')(\log r)^R(\log r')^{R'}(\log h)^H(\log h')^{H'}F_K^s\left (r[h,N],r'[h',N']\right )}{hh'[r[h,N],r'[h',N']]^s}%\hspace {10mm}
%\\ &=:&\sum _{h,h'=1}^\infty \frac {\mu (h)\mu (h')(\log h)^H(\log h')^{H'}\mathcal G_{[\mathbf h,\mathbf N]}(s)}{hh'}\notag 
\end {eqnarray}
%and define so that
%\begin {eqnarray*}
%\frac {d}{dw^Rdw'^{R'}}R_{N,N'}^s(w,w')\vline _{\mathbf w=\mathbf 0}=\mathcal G_\mathbf N(s)
%\end {eqnarray*}
From Lemma \ref {arch} (A) and \eqref {derivative}
\begin {eqnarray*}
\mathcal R_{N,N'}(s)&=&\frac {P^{2K+1,R+R'}(s)F_{K+1}^s\left (N,N'\right )}{[N,N']^s}
\end {eqnarray*}
so 
\begin {eqnarray*}
\mathcal F_\mathbf N(s)&=&P^{2K+1,R+R'}(s)\sum _{h,h'=1}^\infty \frac {\mu (h)\mu (h')(\log h)^H(\log h')^{H'}F_{K+1}^s\left ([h,N],[h',N']\right )}{hh'[[h,N],[h',N']]^s}
%\\ &=:&P_{2K+1}^{R+R'}(s)\mathcal I_\mathbf N(s).
\end {eqnarray*}
and reasoning similarly but using Lemma \ref {arch} (B) allows us to conclude
%Now defining $\mathcal H_{N,N'}^s(w,w')$ as in Lemma \ref {arch} and reasoning similarly, using Lemma \ref {arch} (B), allows us to conclude %\begin {eqnarray*}
%\frac {d}{dw^Hdw'^{H'}}H_{N,N'}^s(w,w')\vline _{\mathbf w=\mathbf 0}%=\mathcal I_\mathbf N(s)
%\end {eqnarray*}
%and we conclude from Lemma \ref {arch} (B)%
%\begin {eqnarray*}
%\mathcal I_\mathbf N(s)&=&\frac {\lambda (s)}{[N,N']^s}F_{K+1}^s\left (\frac {\mathbf N}{\mathcal N,\mathcal N}\right )
%\end {eqnarray*}
%and 
\begin {eqnarray}\label {blodyn}
\mathcal F_\mathbf N(s)&=&\frac {P^{2K+1,R+R'}(s)F_{K+1}^s\left (N,N'\right )}{[N,N']^s}.
\end {eqnarray}
For $n,V,v_1,...,v_V\in \mathbb N\cup \{ 0\} $ write $P_s^{n}(\mathbf v)$ for a quantity of the form
\[ \sum _{k_i,V_i\geq 0\atop {k_1+\cdot \cdot \cdot +k_n+V_1+\cdot \cdot \cdot +V_V\leq d}}\lambda _{\mathbf k,\mathbf V,d,n}(s)\zeta ^{(k_1)}(s)\cdot \cdot \cdot \zeta (s)^{(k_n)}(\log v_1)^{V_1}\cdot \cdot \cdot (\log v_V)^{V_V}\]
so that \eqref {blodyn} says
\begin {eqnarray*}
%\mathcal F_\mathbf N^*(s)&:=&
&&\sum _{r,r',h,h'=1}^\infty \frac {\mu (h)\mu (h')P_s^{n}(\mathbf v,r,r',h,h')F_K^s\left (r[h,N],r'[h',N']\right )}{hh'[r[h,N],r'[h',N']]^s}
%\\ \arrayfi =\sum _{k_i,V_i,R,R',H,H'\geq 0\atop {k_0+\cdot \cdot \cdot +k_d\leq M\atop {\max _{k_j\not =0}j+\mathbf V+\mathbf R+\mathbf H\leq d}}}\lambda _{\mathbf k,\mathbf V,\mathbf R,\mathbf H}(s)\zeta (s)^{\mathbf k}(\log \mathbf v)^{\mathbf V}
%\\ &&\hspace {25mm}\sum _{r,r',h,h'=1}^\infty \frac {\mu (h)\mu (h')(\log \mathbf r)^\mathbf R(\log \mathbf h)^\mathbf HF_K^s\left (r[h,N],r'[h',N']\right )}{hh'[r[h,N],r'[h',N']]^s}
%\\ \arrayfi =\frac {F_{K+1}^s\left (N,N'\right )}{[N,N']^s}\sum _{k_i,V_i,R,R',H,H'\geq 0\atop {k_0+\cdot \cdot \cdot +k_d\leq M\atop {\max _{k_j\not =0}j+\mathbf V+\mathbf R+\mathbf H\leq d}}}\lambda _{\mathbf k,\mathbf V,\mathbf R,\mathbf H}(s)\zeta (s)^{\mathbf k}(\log \mathbf v)^{\mathbf V}P^{2K+1,\max \{ R,R'\} }(s)
%\\ \arrayfi =
=\frac {P_s^{n+2K+1}(\mathbf v)F_{K+1}^s\left (N,N'\right )}{[N,N']^s}.
\end {eqnarray*}
The polynomial in \eqref {spark} is $P_s^{0}(\mathbf r,\mathbf r',\mathbf h,\mathbf h')$ so applying the above equality $k$ times the $r_1,...,h+k'$ sum in \eqref {spark} becomes $P^{k(k+2),d}(s)$ %so this says that the $r_1,...,h_k'$ sum in \eqref {spark} is
%\begin {eqnarray*}
%\sum _{r_1|\cdot \cdot \cdot |r_{k-1}\atop {r_1'|\cdot \cdot \cdot |r_{k-1}'\atop {h_1|r_1,...,h_{k-1}|r_{k-1}\atop {h_1'|r_1',...,h_{k-1}'|r_{k-1}'}}}}\frac {\mu (h_1)\cdot \cdot \cdot \mu (h_{k-1}')F_2^s(r_{k-1},r_{k-1}')}{h_1\cdot \cdot \cdot h_{k-1}'[r_{k-1},r_{k-1}']^s}P_s^{3,...}(\mathbf r,\mathbf r',\mathbf h,\mathbf h')
%\end {eqnarray*}
%and now we can do this another $k-1$ times so that the sum becomes $P^{(k+1)^2-1,d}(s)$ %so that the whole expression in \eqref {spark} is
%\begin {eqnarray*}
%&&\sum _{0\leq d\leq D}\lambda _d(s)\zeta ^{(d)}(s)\sum _{k_i\geq 0\atop {k_0+\cdot \cdot \cdot +k_d\leq (k+1)^2-1}}\lambda _{\mathbf k}(s)\zeta (s)^{k_0}\cdot \cdot \cdot \zeta ^{(d)}(s)^{k_d}
%\end {eqnarray*}
and the result follows.
\end {proof}

\begin {center}
\section {-\hspace {5mm}{Proof of Theorem \ref {mh}}}
\end {center}
Let $V(x,Q)$ be as given in our Theorem \ref {mh} and in view of Theorem \ref {largesieve} take large real numbers $\sqrt x\leq Q=o(x)$.  Let $\gamma $ be a parameter at our disposal subject to $2\sqrt x\leq \gamma \leq x$.  Throughout we will drop $x^\epsilon $ and $\log x$ factors from our error terms and the implied constants will depend on $k,\epsilon $.
%\begin {equation}\label {parameters}
%2\sqrt x\leq \gamma \leq x.
%\end {equation}
\\
\\ On opening up the square in $V(x,Q)$ we have
\begin {eqnarray}\label {variance}
V(x,Q)&=&\sum _{q\leq Q}\sum _{n,m\leq x\atop {n\equiv m(q)}}d_k(n)d_k(m)
%\\ &&-\hspace {4mm}
-2x\sum _{q\leq Q}\frac {1}{q}\sum _{n\leq x}d_k(n)f_x(q,n)+x^2\sum _{q\leq Q}\frac {1}{q^2}\sum _{a=1}^qf_x(q,a)^2\hspace {10mm}
%\\ &=:&S_1(x,Q)-S_2(x,Q).
\end {eqnarray}
where $f_x(q,a)$ is as in Lemma \ref {kfold}.  If $c_q(n)$ is Ramanujan's sum then for $d,d'|q$ 
\begin {eqnarray*}
\sum _{a=1}^qc_d(a)c_{d'}(a)=q\left \{ \begin {array}{ll}\phi (d)&\text { if }d=d'\\ 0&\text { if }d\not =d'\end {array}\right .
\end {eqnarray*}
so 
%\begin {eqnarray*}
%\sum _{a=1}^q\mathcal M_x(q,a)^2&=&q\sum _{d|q}\frac {f_x(d)^2}{d}%=:\mathcal M_x(q)
%\end {eqnarray*} 
\begin {eqnarray*}
\sum _{a=1}^qf_x(q,a)^2&=&q\sum _{d|q}\frac {\phi (d)\mathcal M_x(d)^2}{d^2}=:qf_x^*(q)
\end {eqnarray*}
where, for some polynomial $P$ of degree $\leq 2k-1$,
\begin {eqnarray*}
\sum _{q\leq Q}\frac {f_x^*(q)}{q}=\underbrace {\sum _{d=1}^\infty \frac {\phi (d)f_x(d)^2(\log (Q/d)+\gamma )}{d^3}}_{=P(\log x,\log Q)}+\mathcal O\left (\frac {1}{Q}\right )
\end {eqnarray*}
whilst from Lemma \ref {kfold} (D)
\begin {eqnarray*}
\sum _{n\leq x}d_k(n)f_x(q,n)&=&xf_x^*(q)+\mathcal O\left (x^{1-1/k}\right )
\end {eqnarray*}
We put these and \eqref {square} in \eqref {variance} to get
%\[ \sum _{q\leq Q}\sum _{n,m\leq x\atop {n=m}}d_k(n)d_k(m)=xQP(\log x)+\mathcal O\left (\right )\]
\begin {eqnarray}\label {initialdispositions}
V(x,Q)&=&2\sum _{q\leq Q}\sum _{m<n\leq x\atop {n\equiv m(q)}}d_k(n)d_k(m)-xQP(\log x)+x^2P(\log x,\log Q)+\mathcal O\left (x^{2-1/k}+Qx^{1-\mathfrak d}\right ).\hspace {10mm}
\end {eqnarray}
As said in the introduction, we now borrow from \cite {goldstonvaughan}.  There are, for each $1\leq a\leq q\leq \gamma $ with $(a,q)=1$, disjoint intervals $\mathfrak F(a/q)$ about $a/q$ such that for any continuous $f:\mathbb R\rightarrow \mathbb C$ of period 1
\[ \int _0^1f(t)dt=\sum _{q\leq \gamma }\sideset {}{'}\sum _{a=1}^q\int _{\mathfrak F(a/q)}f(t)dt.\]
The intervals satisfy
\begin {equation}\label {farey}
\left (\frac {a}{q}-\frac {1}{2q\gamma },\frac {a}{q}+\frac {1}{2q\gamma }\right )\subseteq \mathfrak F(a/q)\subseteq \left (\frac {a}{q}-\frac {1}{q\gamma },\frac {a}{q}+\frac {1}{q\gamma }\right ).
\end {equation}
This is the \emph {Farey dissection} of the unit interval of order $\gamma$ and for its discussion see, for example, Section 3.8 of \cite {hardywright}.  Defining for $X>0$, $\alpha \in \mathbb R$
\[ S_X(\alpha )=\sum _{n\leq X}d_k(n)e(n\alpha )\hspace {10mm}\text {and}\hspace {10mm}F(\alpha )=\sum _{uv\leq x\atop {u\leq Q}}e\left (uv\alpha \right )\]
we see that
\begin {eqnarray}\label {circlemethod}
\sum _{uv\leq x\atop {u\leq Q}}\sum _{m<n\leq x\atop {n-m=uv}}d_k(n)d_k(m)=\sum _{q\leq \gamma }\sideset {}{'}\sum _{a=1}^q\int _{\mathfrak F(a/q)}F(\alpha )|S_x(\alpha )|^2d\alpha =:\mathfrak I
\end {eqnarray}
and as the LHS here is the $q,n,m$ sum in \eqref {initialdispositions} we can say
\begin {equation}\label {hwn}
V(x,Q)=2\mathfrak I+xQ\underbrace {P(\log x)}_{\text {degree }k^2-1}+x^2\underbrace {P(\log x,\log Q)}_{\text {degree }\leq 2k-1}+\mathcal O\left (x^{2-1/k}+Qx^{1-\mathfrak d}\right ).
\end {equation}
Define
\[ F_q(\alpha )=\sum _{u\leq \sqrt x\atop {q|u}}\left (\sum _{v\leq x/u}+\sum _{\sqrt x<v\leq Q,x/u}\right )e(uv\alpha )\]
and write $\alpha =a/q+\beta $.  From (3.1), (3.4), (3.5) and (3.6) of \cite {goldstonvaughan} and from \eqref {farey} we have $F(\alpha )=F_q(\alpha )+\mathcal O\left (\gamma \right )$ for $\alpha \in \mathfrak F(a/q)$ and from (3.11) of \cite {goldstonvaughan} also 
\begin {eqnarray}\label {FF}
F_q(\alpha )&\ll &\frac {x}{q(1+x|\beta |)}\hspace {10mm}|\beta |\leq \frac {1}{2q\sqrt x}
\\ &\ll &\gamma \hspace {10mm}\text { for }\alpha \in \mathfrak F(a/q)\char92 \left (\frac {a}{q}-\frac {1}{2q\gamma },\frac {a}{q}+\frac {1}{2q\gamma }\right )\notag 
\end {eqnarray}
so we can say\footnote {here and on note that $F_q(\beta )=0$ once $q>\sqrt x$} 
%\hspace {10mm}\text {for }\hspace {2mm}\alpha \in \mathfrak M(a/q)\hspace {3mm}%q\leq \gamma 
\begin {eqnarray}\label {iii}
\mathfrak I&=&\sum _{q=1}^\infty \sideset {}{'}\sum _{a=1}^q\int _{\pm 1/2q\gamma }F_q(\beta )|S_x(\alpha )|^2d\beta +\mathcal O\left (x\gamma \right )
\end {eqnarray}
%Note 
%\begin {equation}\label {ortho}
%\sum _{q\leq \gamma }\sideset {}{'}\sum _{a=1}^q\int _{\mathfrak M(a/q)}|S_x(\alpha )|^2d\alpha \ll x.
%\end {equation}
%For $\alpha \in \mathfrak M(a/q)\char92 (a/q-1/2q\gamma ,a/q+1/2q\gamma )$ we have
%\begin {eqnarray*}
%F(\alpha )&\ll &\gamma 
%\end {eqnarray*}
%so
%Write $\alpha =a/q+\beta $.  For $u\leq \sqrt x$ we have $|u\beta |\leq 1/2q$, therefore $||\alpha u||\gg ||au/q||$, so
%\begin {eqnarray}\label {FF}
%F(\alpha )&=&\sum _{u\leq \sqrt x\atop {q|u}}\left (\sum _{v\leq x/u}+\sum _{\sqrt x<v\leq Q,x/u}\right )e(uv\alpha )+\sum _{q\nmid u}\cdot \cdot \cdot \notag 
%\\ &=:&\underbrace {F_q(\beta )}_{\ll 1/q|\beta |}+\mathcal O\left (\sqrt x+q\right )
%\end {eqnarray}
%so in $\mathfrak I$ we can, using \eqref {ortho}, replace $F(\alpha )$ with $F_q(\beta )$ at the cost of a total error $\ll x^{3/2}+x\gamma $.  Using \eqref {ortho} again and using the bound in \eqref {FF}, we also see that the total contribution to $\mathfrak I$ of the parts $|\beta |\geq 1/qw$ is $\ll xw$.  So in total
%\begin {eqnarray}\label {iii}
%\mathfrak I&=&\sum _{q\leq \gamma }\sideset {}{'}\sum _{a=1}^q\int _{\pm 1/q\gamma }F(\beta )|S_x(\alpha )|^2d\beta +\mathcal O\left (x\gamma \right )
%\end {eqnarray}
and now we do some book-keeping before turning to this main term.  Define%Take $P(\cdot ,...,\cdot )$ to be the polynomial of Lemma \ref {kfold} (B) and define
\begin {eqnarray*}
%\text {for }X>0,\hspace {1mm}q\in \mathbb N\hspace {10mm}
\mathcal M_X(q)&&\text {as in Lemma \ref {kfold}}\hspace {15mm}
g_q(t)=\frac {d}{dt}\{ t\mathcal M_t(q)\} \hspace {15mm}\\ 
%\\ \text {for }X>0,\hspace {1mm}a,q\in \mathbb N\hspace {10mm}
%S_X(a/q)-\frac {X\mathcal M_X(q)}{q}\notag 
%\\ \text {for }a,q\in \mathbb N\hspace {10mm}\hat \Delta (a/q)&=&\max _{X\leq x}|\Delta _X^k(a/q)|\notag 
%\\ \text {for }q\in \mathbb N,\hspace {1mm}t\in \mathbb R\hspace {10mm}
\Delta _X(a/q)&&\text {as in Lemma \ref {barbandavenporthalberstam}}\notag \hspace {15mm}
%\\ \text {for }q\in \mathbb N,\hspace {1mm}\beta \in \mathbb R\hspace {10mm}
I_q(\beta )=\int _1^xe(\beta t)g_q(t)dt
\end {eqnarray*}
and write for $k>2$
\[ \sum _{\mathbf h|\mathbf r|q\atop {\mathbf h'|\mathbf r'|q}}=\sum _{r_1|\cdot \cdot \cdot |r_{k-2}|q\atop {r_1'|\cdot \cdot \cdot |r_{k-2}'|q\atop {h_1|r_1,...,h_{k-2}|r_{k-2}\atop {h_1'|r_1',...,h_{k-2}'|r_{k-2}'}}}}\frac {\mu (h_1)\cdot \cdot \cdot \mu (h'_{k-2})}{h_1\cdot \cdot \cdot h_{k-2}'}\]
and for $k=2$ 
\[ \sum _{\mathbf h|\mathbf r|q\atop {\mathbf h'|\mathbf r'|q}}P_{D,D'}(q,\mathbf r,\mathbf r',\mathbf h,\mathbf h')=P_{D,D'}(q).\]
Note that there are polynomials $P_{D,D'}(\cdot ,...,\cdot )$ in $\log (\cdot ),...,\log (\cdot )$ of degree $\leq 2k-2-D-D'$ such that
\begin {eqnarray}\label {cardiadolig}
g_q(t)g_q(t')&=&\sum _{0\leq D,D'\leq k-1}(\log t)^{D}(\log t')^{D'}\sum _{\mathbf h|\mathbf r|q\atop {\mathbf h'|\mathbf r'|q}}P_{D,D'}(q,\mathbf r,\mathbf r',\mathbf h,\mathbf h')\notag 
\\ &=:&\sum _{0\leq D,D'\leq k-1}(\log t)^{D}(\log t')^{D'}\alpha _{\mathbf D}(q)
\end {eqnarray}
and that for $t\leq x$
\begin {eqnarray}\label {ibound}
I_q(\beta )\ll \frac {x}{1+x|\beta |}\hspace {8mm}g_q(t)\ll 1\hspace {8mm}g_q'(t)\ll \frac {1}{t}.
\end {eqnarray}
From Lemma \ref {barbandavenporthalberstam} and the Cauchy-Schwarz inequality  %Note also that for $k=2$ Lemma \ref {kfold} (B) says  %$\hat \Delta _x(q)\ll x^{1-2/(2k-1)}$
%\begin {equation}\label {k2} %\sum _{q\leq X}\frac {1}{q}\sum _{a=1}^q|\Delta _y^2(a/q)|^2\ll Xy^{2/3}+X^2\hspace {5mm}\text { and }\hspace {5mm}
%\sum _{q\leq \gamma }\frac {1}{q^2}\sum _{a=1}^q|\Delta _x^2(a/q)|^2\ll x
%\end {equation}
%Note also that for $k=2$ Lemma \ref {kfold} (B) says for $0<X\leq \sqrt y$ %$\hat \Delta _x(q)\ll x^{1-2/(2k-1)}$
%\begin {equation}\label {k2} %\sum _{q\leq X}\frac {1}{q}\sum _{a=1}^q|\Delta _y^2(a/q)|^2\ll Xy^{2/3}+X^2\hspace {5mm}\text { and }\hspace {5mm}
%\sum _{q\leq X}\frac {1}{q^2}\sum _{a=1}^q|\Delta _y^2(a/q)|^2\ll y.
%\end {equation}
%and, defining $E_x(q,a)$ and $\theta _k$ as in our theorem,
\begin {eqnarray}\label {arcs2}
\sum _{q\leq \sqrt x}\frac {1}{q^2}\sideset {}{'}\sum _{a=1}^q\max _{t\leq x}|\Delta _t(a/q)|\ll \sqrt x+\underbrace {x^{1-4/(6k-3)}}_{k>2}
\end {eqnarray}
and %from Lemma \ref {barbandavenporthalberstam}
\begin {eqnarray}\label {arcs}
\sum _{q\leq \sqrt x}\frac {1}{q}\left (1+\frac {x}{q\gamma }\right )\sideset {}{'}\sum _{a=1}^q\max _{t\leq x}|\Delta _t(a/q)|^2\ll x^{3/2}+\underbrace {x^{2-1/k}+\frac {x^{3-8/(6k-3)}}{\gamma }}_{k>2}.
\end {eqnarray}
%so \textbf {if our theorem is true for some $k$} then for $0<X\leq y$
%\begin {eqnarray}\label {kgeneral}
%\sum _{q\leq X}\frac {1}{q}\sum _{a=1}^q|\Delta _y^k(a/q)|^2\ll Xy+y^{2-\Delta _{k}}\hspace {5mm}\text { and }\hspace {5mm}\sum _{q\leq X}\frac {1}{q^2}\sum _{a=1}^q|\Delta _y^k(a/q)|^2\ll y^{2-\Delta _{k}}\hspace {10mm}
%\end {eqnarray}
%where $\theta _k$ are as given in our theorem.  
Now back to \eqref {iii}.  Partial summation, integration by parts and \eqref {ibound} give
\begin {eqnarray*}
S_x(\alpha )&=&e(x\beta )S_x(a/q)-2\pi i\beta \int _1^xe(t\beta )S_t(a/q)dt
\\ &=&\underbrace {\frac {I_q(\beta )}{q}}_{\ll x/q(1+|\beta |x)}+\underbrace {e(x\beta )\Delta _x(a/q)-2\pi i\beta \int _1^xe(t\beta)\Delta _t(a/q)dt}_{=:J(\alpha )\ll (1+|\beta |x)\max _{t\leq x}|\Delta _t(a/q)|}+\mathcal O(1)
\end {eqnarray*}
so
\begin {eqnarray}\label {sss}
|S_x(\alpha )|^2&=&\frac {|I_q(\beta )|^2}{q^2}+\mathcal O\left (\frac {x\max _{t\leq x}|\Delta _t(a/q)|}{q}+|J(\alpha )|^2+1\right ).
\end {eqnarray}
%and if $k=2$ we can use [jutila] in the last but one equality and say the error term is actually
%\begin {eqnarray*}
%&\ll &\frac {x^{5/2}}{q\gamma ^2}\hspace {10mm}\text { if $\gamma \ll \sqrt x$}.
%\end {eqnarray*}
With \eqref {FF} we have
\begin {eqnarray*}
\int _{\pm 1/2q\gamma }F_q(\beta )|J(\alpha )|^2d\beta 
&\ll &
%\\ \arrayfi \ll 
\frac {|\Delta _x(a/q)|^2}{q}\int _{\pm 1/q\gamma }\frac {xd\beta }{1+x|\beta |}+\frac {1}{q}\int _1^x\int _1^x\Delta _t(a/q)\overline {\Delta _{t'}(a/q)}\int _{\pm 1/q\gamma }\beta e\left (\beta (t-t')\right )dt'dt
%\\ \arrayfi \ll 
\\ &\ll &\frac {1}{q}\left (1+\frac {x}{q\gamma }\right )\max _{t\leq x}|\Delta _t(a/q)|^2
\end {eqnarray*}
so from \eqref {sss}, \eqref {FF}, \eqref {arcs2} and \eqref {arcs} the main term in \eqref {iii} is
\begin {eqnarray*}
%&&\sum _{q\leq \gamma }\sideset {}{'}\sum _{a=1}^q\int _{\pm 1/qw}F_q(\alpha )|S_x(\alpha )|^2d\alpha 
%\\ \arrayfi =
\sum _{q=1}^\infty \frac {\phi (q)}{q^2}\int _{\pm 1/2q\gamma }F_q(\beta )|I_q(\beta )|^2d\beta +\mathcal O\left (x^{3/2}+\underbrace {x^{2-4/(6k-3)}+\frac {x^{3-8/(6k-3)}}{\gamma }}_{k>2}\right ).
%x^{2-1/k}+\sum _{q\leq \gamma }\frac {1}{q}\sum _{a=1}^q\max _{X\leq x}|\Delta _X^k(a/q)|^2\left (1+\frac {x}{qw}\right) \right ).
%\\ \arrayfi {=:}M_x(\gamma )+\mathcal O\left (\mathfrak E\right ).
\end {eqnarray*}
From \eqref {ibound} and $F_q(\beta )\ll x/q$ the integral here may be extended to infinity at the cost of an error $\ll x\gamma $ and we conclude
%For any $X_{q}\subseteq [0,1]$ and any $Z>0$ the bounds in \eqref {ibound} and \eqref {FF} say
%\begin {eqnarray*}
%\sum _{q>Z}\frac {1}{q}\int _{X_q}F_q(-\beta )|I_q(\beta )|^2d\alpha %&\leq &x^2\sum _{q>Z}\frac {1}{q^2}
%\\ &\ll &
%\ll \frac {x^2}{Z}
%\end {eqnarray*}
%so from \eqref {ibound}
%\begin {eqnarray*}
%\sum _{q\leq \gamma }\frac {1}{q}\int _{1/q\gamma }^{1/2}F_q(-\beta )|I_q(\beta )|^2d\alpha &\ll &%\sum _{q\leq x/w}\frac {1}{q^2}\int _{1/q\gamma }^\infty \frac {d\beta }{\beta ^3}+x\gamma 
%\\ &\ll &
%\ll x\gamma 
%\end {eqnarray*}
%so the integral above is
%\begin {eqnarray*}
%\sum _{q\leq \gamma }\int _{\pm 1/qw}%=\sum _{q=1}^\infty \int _{\pm 1/2}+\mathcal O\left (\sum _{q\leq \gamma }\int _{1/qw}^{1/2}+\sum _{q>\gamma }\int _{\pm 1/2}\right )
%=\sum _{q=1}^\infty \int _{\pm 1/2}+\mathcal O\left (x\gamma +\frac {x^2}{\gamma }\right )
%\end {eqnarray*}
%and %setting $\gamma =\sqrt x$ for $k=2$ and $\gamma =x^{1-\theta _k/2}$ for $k>2$
\begin {eqnarray}\label {eee}
\mathfrak I&=&\underbrace {\sum _{q=1}^\infty \frac {\phi (q)}{q^2}\int _{\pm \infty }F_q(\beta )|I_q(\beta )|^2d\beta }_{=:\mathfrak M}+\mathcal O\left (x^{3/2}+\underbrace {x^{2-4/(6k-3)}}_{k>2}\right ).%\notag 
%\\ &=:&\mathfrak M+\mathcal O\left (x^{3/2}+\underbrace {x^{2-1/2k}}_{\text {if }k>2}\right ).
\end {eqnarray}
%\begin {eqnarray}\label {eee}
%\mathfrak I&=&\sum _{q=1}^\infty \frac {\phi (q)}{q^2}\int _{\pm 1/2}F_q(\beta )|I_q(\beta )|^2d\beta +\notag 
%\\ \arrayfi =\mathcal O\left (xw+\frac {x^2}{\gamma }+x^{2-1/k}+\sum _{q\leq \gamma }\frac {1}{q}\sum _{a=1}^q\max _{X\leq x}\left |\Delta _X^k(a/q)\right |^2\left (1+\frac {x}{qw}\right )\right )\notag 
%\\ &=:&\mathfrak M+\mathcal O\left (xw+\frac {x^2}{\gamma }+x^{2-1/k}+\mathfrak E_k\right ).
%\end {eqnarray}
Let 
\begin {eqnarray*}
L(Z)&=&\int _{1,1+Z}^{x,x+Z}g_q(t)g_q(t-Z)dt
%\\ &=&\left \{ \begin {array}{ll}&\text { for \end {array}\right .
\end {eqnarray*}
so that $|I_q(\beta )|^2=\hat L(\beta )$ and therefore
\begin {eqnarray*}
\int _{\pm \infty }e(-uv\beta )|I_q(\beta )|^2d\beta =L(uv)
\end {eqnarray*}
%making the integral in $\mathfrak M$
%\begin {eqnarray*}
%&&\sum _{u\leq \sqrt x\atop {q|u}}\left (\sum _{v\leq x/u}+\sum _{\sqrt x<v\leq Q,x/u}\right )\int _{\pm \infty }e(-uv\beta )|I_q(\alpha )|^2d\beta 
%\\ \arrayfi =\sum _{u\leq \sqrt x\atop {q|u}}\left (\sum _{v\leq x/u}+\sum _{\sqrt x<v\leq Q,x/u}\right )\int _1^x\int _1^xg_q(t)g_q(t')\left (\int _{\pm \infty }e\left (\beta (t-t'-uv)\right )d\beta \right )dt'dt
%\\ \arrayfi =\sum _{u\leq \sqrt x\atop {q|u}}\left (\sum _{v\leq x/u}+\sum _{\sqrt x<v\leq Q,x/u}\right )\int _{1+uv}^xg_q(t)g_q(t-uv)dt
%\end {eqnarray*}
so
\begin {eqnarray}\label {cowbois}
\mathfrak M&=&\sum _{q=1}^\infty \frac {\phi (q)}{q^2}\sum _{u\leq \sqrt x\atop {q|u}}\int _{1+u}^xg_q(t)\left (\sum _{v\leq (t-1)/u}+\sum _{\sqrt x<v\leq Q,(t-1)/u}\right )g_q(t-uv)dt\notag 
\\ &=:&\sum _{q=1}^\infty \frac {\phi (q)}{q^2}\sum _{u\leq \sqrt x\atop {q|u}}\int _{1+u}^xg_q(t)V_u(t).\hspace {10mm}
%\\ &=:&\sum _{q=1}^\infty \frac {S_q(x,Q)}{q^2}+\mathcal O\left (\mathfrak E_k\right ).
\end {eqnarray}
Write
\[ G_q(X)=\int _1^Xg_q(t)dt%\hspace {6mm}\text { and }\hspace {6mm}V_u(t)=\left (\sum _{v\leq (t-1)/u}+\sum _{\sqrt x<v\leq Q,(t-1)/u}\right )g_q(t-uv)
\]
and bear \eqref {ibound} in mind.  For any $0\leq A<B\leq (t-1)/u$ the Euler-Maclaurin summation formula says that
\begin {eqnarray*}
\sum _{A<v\leq B}g_q(t-uv)=\frac {G_q(t-uA)-G_q(t-uB)}{u}+\mathcal O(1)
\end {eqnarray*}
so the first sum in $V_u(t)$ is
\[ G_q(t)/u+\mathcal O(1)\]
and the second sum vanishes for $(t-1)/u<\sqrt x$ whilst for $(t-1)/u\geq \sqrt x$ it's
\begin {eqnarray*}
%&&-\frac {1}{u}\left \{ \begin {array}{ll}F(t)&\text { if }t/u\leq Q\\ F(t)-F(t-uQ)&\text { if }t/u>Q\end {array}\right \} +F\left (t-u\sqrt x\right )+\mathcal O\left (1\right )
%\\ &&\hspace {10mm}=\hspace {4mm}
\frac {G_q\left (t-u\sqrt x\right )}{u}-\frac {1}{u}\left \{ \begin {array}{ll}0&\text { if }(t-1)/u\leq Q\\ G_q(t-uQ)&\text { if }(t-1)/u>Q\end {array}\right \} +\mathcal O\left (1\right )
\end {eqnarray*}
so \eqref {cowbois} becomes
\begin {eqnarray}\label {cowbois2}
\mathfrak M&=&\sum _{q=1}^\infty \frac {\phi (q)}{q^3}\left (\sum _{u\leq \sqrt x/q}\frac {1}{u}\right )\int _1^xg_q(t)G_q(t)dt+\{ \text { term $Q=\sqrt x$ }\} \notag 
\\ &&-\hspace {4mm}\int _{1+Q}^xg_q(t)\left (\sum _{u\leq (t-1)/Q}\frac {1}{u}\sum _{q|u}\frac {\phi (q)}{q^2}G_q(t-uQ)\right )dt+\mathcal O\left (x^{3/2}\right )\notag 
\\ &=:&\mathcal A(x)+\mathcal B(\sqrt x)-\mathcal B(Q)+\mathcal O\left (x^{3/2}\right ).
\end {eqnarray}
From \eqref {cardiadolig} the integrand in $\mathcal B(Q)$ is 
\begin {eqnarray*}
&&\sum _{0\leq D,D'\leq k-1}(\log t)^{D}\sum _{u\leq (t-1)/Q}\frac {1}{u}\underbrace {\left (\sum _{q|u}\frac {\phi (q)\alpha _\mathbf D(q)}{q^2}\right )}_{=:\Phi _\mathbf D(u)}\int _1^{t-uQ}(\log t')^{D'}dt'
%\\ \arrayfi {=:}\sum _{0\leq D,D'\leq k-1}(\log t)^{D}\sum _{u\leq (t-1)/Q}\frac {\Phi _{\mathbf D}(u)}{u}\int _1^{t-uQ}(\log t')^{D'}dt'
\end {eqnarray*}
so from Lemma \ref {w} for any $c>0$
\begin {eqnarray}\label {bq}
\mathcal B(Q)&=&\sum _{0\leq D,D'\leq k-1}\int _{c\pm \infty }\frac {\mathcal F_\mathbf D(s)}{Q^s}\left (\int _1^x(\log t)^{D}K_t(s)dt\right )ds+\mathcal O\left (Q^2\right )\hspace {10mm}
\end {eqnarray}
where
\begin {eqnarray*}
\mathcal F_\mathbf D(s)&=&\sum _{u=1}^\infty \frac {\Phi _{\mathbf D}(u)}{u^{s+1}}
\end {eqnarray*}
and where $K_t(s)$ satisfies
\begin {eqnarray*}
\int _1^x(\log t)^{D}K_t(s)dt&\ll &\frac {x^{\sigma +2}}{|s(s+1)(s+2)|}
\end {eqnarray*}
amd is holomorphic for $\sigma >-2$ except for simple poles at $s=-1,0$ where it has Laurent expansions
\begin {eqnarray*}
\frac {K_t(s)}{Q^s}%&=&
&=&\frac {t}{s}\sum _{n\geq 0}P_n(t,Q)s^n\hspace {10mm}\text {about }s=0
%\\ 
%&=&\frac {x^2+\mathcal O\left (Q^2\right )}{s}\sum _{n\geq 0}P_n(\log x,\log Q)s^n\hspace {10mm}\text {about }s=0
\\ \frac {K_t(s)}{Q^s}%&=&
&=&\frac {Q}{s+1}\sum _{n\geq 0}P_n(t,Q)(s+1)^n\hspace {10mm}\text {about }s=-1
%&=&\frac {xQ+\mathcal O\left (Q^2\right )}{s+1}\sum _{n\geq 0}P_n(\log x,\log Q)(s+1)^n\hspace {10mm}\text {about }s=-1.
\end {eqnarray*}
for some polynomials in $\log t,\log Q$ of degree $D'+n$.  Recall the definitions in the sentence containing \eqref {cardiadolig} and write $Z=k(k-2)$.  From Lemma \ref {arch2} there are $\lambda _\mathbf k(s)$ holomorphic and $\ll 1$ for $\sigma >-3/2$ such that for $s>-1$
\begin {eqnarray*}
\mathcal F_\mathbf D(s)&=&\zeta (s+1)\sum _{q=1}^\infty \frac {\phi (q)}{q^{s+3}}\sum _{r_1|\cdot \cdot \cdot |r_{k-2}|q\atop {r_1'|\cdot \cdot \cdot |r_{k-2}'|q\atop {h_1|r_1,...,h_{k-2}|r_{k-2}\atop {h_1'|r_1',...,h_{k-2}'|r_{k-2}'}}}}\frac {\mu (h_1)\cdot \cdot \cdot \mu (h_{k-2}')}{h_1\cdot \cdot \cdot h_{k-2}'}\underbrace {P_\mathbf D(q,\mathbf r,\mathbf r',\mathbf h,\mathbf h')}_{\text {degree }\leq 2k-2-D-D'=:d}
\\ &=&\zeta (s+1)\sum _{k_i\geq 0\atop {k_0+\cdot \cdot \cdot +k_Z\leq d}}\lambda _{\mathbf k}(s)\zeta ^{(k_1)}(s+2)\cdot \cdot \cdot \zeta ^{(k_Z)}(s+2)
%\\ &=&\zeta (s+1)\sum _{k_i\geq 0\atop {k_0+\cdot \cdot \cdot +k_{d}\leq k(k-2)}}\lambda _\mathbf k(s)\zeta (s+2)^{k_0}\cdot \cdot \cdot \zeta ^{(d)}(s+2)^{k_{d}}
\\ &=:&\zeta (s+1)\lambda (s)
\\ &\ll _k&|\zeta (s+1)|\max _{0\leq n\leq d}|\zeta ^{(n)}(s+2)|^{Z}
\end {eqnarray*}
so the integral in \eqref {bq} may be moved to $\mathfrak c-2$, where $\mathfrak c$ is as in Theorem \ref {mh}, where it is bounded by
\begin {eqnarray*}
%\int _{\mathfrak c_{n,d}-2}\frac {|\zeta (s+1)\mathcal F_\mathbf T(s)}{t^2}&\ll &\frac {x^{s+2}}{Q^s}\int _{\mathfrak c_{n,d}-2}\frac {t^{1/2-(\mathfrak c_{n,d}-2+1)}|\zeta ^{(n)}(s+2)|^{d}}{t^3}
%\\ &\ll &
\frac {x^{\mathfrak c}}{Q^{\mathfrak c-2}}\int _{\mathfrak c-2\pm \infty }\frac {|\zeta ^{(n)}(s+2)|^{k(k-2)}ds}{(1+t)^{3/2+\mathfrak c}}
%\\ &\ll &
\ll Q^2\left (\frac {x}{Q}\right )^\mathfrak c,
\end {eqnarray*}
picking up residues
\begin {eqnarray*}
\Big (Res _{s=0}+Res _{s=-1}\Big )\Bigg \{ \frac {\zeta (s+1)\lambda (s)}{Q^s}\int _1^x(\log t)^{D}K_t(s)dt\Bigg \} .
\end {eqnarray*}
From the above Laurent expansions the first residue is
\begin {eqnarray*}
\Big (x^2+\mathcal O(1)\Big )Res_{s=0}\Bigg \{ \frac {\zeta (s+1)\lambda (s)}{s}\sum _{n\geq 0}s^n\underbrace {P_n(x,Q)}_{\text {degree }\leq n+D+D'}\Bigg \} 
\end {eqnarray*}
and the second is
\begin {eqnarray*}
Q\Big (x+\mathcal O(1)\Big )Res_{s=-1}\Bigg \{ \frac {\zeta (s+1)}{s+1}\left (\frac {1}{(s+1)^{Z+d}}\sum _{k_i\geq 0\atop {k_0+\cdot \cdot \cdot +k_Z=d}}\lambda _\mathbf k(s)+\cdot \cdot \cdot \right )\sum _{n\geq 0}(s+1)^n\underbrace {P_n(x,Q)}_{\text {degree }\leq n+D+D'}\Bigg \} 
\end {eqnarray*}
so the sum of these residues is $x^2P(x,Q)+xQP(x,Q)+\mathcal O\left (Q\right )$ for some polynomials in $\log x,\log Q$ of degrees $\leq 1+D+D'$ and $\leq Z+d+D+D'$ respectively and \eqref {bq} becomes
\[ \mathcal B(Q)=x^2P(x,Q)+xQP(x,Q)+\mathcal O\left (Q^2\left (\frac {x}{Q}\right )^\mathfrak c\right )\]
for some polynomials of degree $\leq 2k-1$ and $\leq k^2-1$.  Meanwhile it is clear that the quantity $\mathcal A(x)$ in \eqref {cowbois2} is $x^2P(x,Q)$ for some polynomial of degree $\leq 2k-1$ so putting all this in \eqref {cowbois2} we get Theorem \ref {mh} from \eqref {hwn} and \eqref {eee}, with Theorem \ref {largesieve} telling us there can be no $x^2$ term. 
\\
\\ 
\begin {center}

\end {center}

\hspace {1mm}
\\
\\
\\
\\ 
\emph {e-mail address} - tomos.parry1729@hotmail.co.uk


\begin {thebibliography}{1}
\bibitem {banks}
W.D. Banks, R. Heath-Brown and I.E. Shparlinski - \emph {On the average value of divisor sums in arithmetic progressions} - International Mathematics Research Notices, 1 (2005)
\bibitem {blomer}
V. Blomer - \emph {The average value of divisor sums in arithmetic progressions} - The Quarterly Journal of Mathematics, 59 (2008)
\bibitem {bf}
R. de la Bretèche , D. Fiorilli - \emph {Major arcs and moments of arithmetical sequences} - American Journal of Mathematics, 142 (2020)
\bibitem {bruedernvaughan}
J. Br\" udern and R. C. Vaughan - \emph {A Montgomery-Hooley theorem for sums of two cubes} - European Journal of
Mathematics (2021)
%\bibitem {bruedern}
%J. Brüdern - \emph {Einführung in die Analytische Zahlentheorie} - Springer, Berlin (1995)
\bibitem {fouvry}
E. Fouvry - \emph {Autour du th\' eor\` eme de Bombieri-Vinogradov} - Acta Mathematica 152 (1984)
%\bibitem {fi}
%J. Friedlander and H. Iwaniec - \emph {The divisor problem in arithmetic progressions} - Acta Arithmetica, XLV (1985)
\bibitem {hb}
R. Heath-Brown - \emph {The divisor function $d_3(n)$ in arithmetic progressions} - Acta Arithmetica, XLVII (1986)
\bibitem {goldstonvaughan}
D. A. Goldston and R. C. Vaughan - \emph {On the Montgomery–Hooley asymptotic formula; in Sieve Methods, Exponential Sums and Their Applications in Number Theory} - Cambridge University Press (1997)
\bibitem {keating}
J. Keating, B. Rodgers, E. Roditty-Gershon and Z. Rudnick - \emph {Sums of divisor functions in $F_q[t]$ and matrix integrals} - Mathematische Zeitschrift, 288 (2018)
\bibitem {lauzhao}
Y.-K. Lau and L. Zhao - \emph {On a variance of Hecke eigenvalues in arithmetic progressions} - Journal of Number Theory, 132 (2012)
\bibitem {hardywright}
G. H. Hardy and E.M. Wright - \emph {The theory of numbers} (3rd edition) - Claredon Press (1954)
%\bibitem {hooley}
%C. Hooley - \emph {On the Barban-Davenport-Halberstam Theorem X} - Hardy-Ramanujan Journal, Volume 21 (1998)
%\bibitem {vaughangeneral}
%\bibitem {hooley3}
%C. Hooley - \emph {On the Barban–Davenport Halberstam theorem III} - Journal of the London Mathematical Society, Volume 10, Issue 1 (1975)
%R. C. Vaughan - \emph {On a variance associated with the distribution of general sequences in arithmetic progressions I} - Philosophical Transactions of the Royal Society of London, Series A (1998)
\bibitem {motohashi}
Y. Motohashi - \emph {On the distribution of the divisor function in arithmetic progressions} - Acta Arithmetica XXII (1973)
\bibitem {nguyen1} 
D. Nguyen - \emph {Variance of the k-fold divisor function in arithmetic progressions for individual modulus} - arXiv
\bibitem {nguyen2} 
D. Nguyen - \emph {Generalized divisor functions in arithmetic progressions: II} -   
\\ https://web.math.ucsb.edu/~dnguyen/preprints/DivisorFunctionsInAPII.pdf
\bibitem {prapanpong}
P. Pongsriiam - \emph {The distribution of the divisor function in arithmetic progressions} - Ph.D. thesis, Pennsylvania State University (2012)
\bibitem {vaughandap}
P. Pongsriiam and R.C. Vaughan - \emph {The divisor function on residue classes I} - Acta Arithmetica, 168 (2015)
\bibitem {vaughandap2}
P. Pongsriiam and R.C. Vaughan - \emph {The divisor function on residue classes II} - Acta Arithmetica, 182 (2018)
\bibitem {sound}
B. Rodgers and K. Soundararajan - \emph {The variance of divisor sums in arithmetic progressions} - Forum Mathematicum 30 (2018)
\bibitem {rzf}
E. C. Titchmarsh - \emph {The theory of the Riemann zeta function} (2nd edition) - Claredon Press, Oxford (1986) 
\bibitem {survey}
R.C. Vaughan - \emph {Generalized Montgomery-Hooley formula; A survey} - 
\\ http://www.personal.psu.edu/rcv4/MontgomeryHooley.pdf
\end {thebibliography}
\end {document}